\documentclass[12pt,a4paper]{article}

\usepackage[margin=1in]{geometry}  
\usepackage{graphicx}              
\usepackage{amsmath,euscript}
\usepackage{amsfonts}              
\usepackage{amssymb,,amscd,epsf,amsthm, epsfig}                
\usepackage[toc,page]{appendix}

\newtheorem{thm}{Theorem}[section]
\newtheorem{fact}[thm]{Fact}
\newtheorem{prop}[thm]{Proposition}
\newtheorem{cor}[thm]{Corollary}
\newtheorem{lemma}[thm]{Lemma}
\newtheorem{claim}[thm]{Claim}
\theoremstyle{definition}
\newtheorem{rem}[thm]{Remark}
\newtheorem{defn}[thm]{Definition}

\input xy
\xyoption{all}

\newcommand{\Z}{\mathbb{Z}}

\newcommand{\calF}{\mathcal{F}}

\newcommand{\calH}{\mathcal{H}}

\newcommand{\DD}{\EuScript D}

\newcommand{\pf}{{\it Proof.}\hspace{2ex}}

\newcommand{\epf}{\hspace*{\fill}\mbox{$\halmos$}}
\newcommand{\halmos}{\rule{1ex}{1.4ex}}

\newcommand{\Hom}{\mathop{\mathrm{Hom}}\nolimits}

\newcommand{\Imm}{\mathop{\mathrm{Im}}\nolimits}
\newcommand{\End}{\mathop{\mathrm{End}}\nolimits}

\newcommand{\Barr}{\mathop{\mathrm{Bar}}\nolimits}
\newcommand{\modu}{\mathop{\mathrm{mod}}\nolimits}

\newcommand{\id}{\mathop{\mathrm{id}}\nolimits}

\newcommand{\op}{\mathop{\mathrm{op}}\nolimits}

\newcommand{\Ker}{\mathop{\mathrm{Ker}}\nolimits}
\newcommand{\Tor}{\mathop{\mathrm{Tor}}\nolimits}
\newcommand{\Tot}{\mathop{\mathrm{Tot}}\nolimits}
\newcommand{\Ext}{\mathop{\mathrm{Ext}}\nolimits}
\newcommand{\sg}{\mathop{\mathrm{sg}}\nolimits}

\newcommand{\HC}{\mathop{\mathrm{HC}}\nolimits}
\newcommand{\HH}{\mathop{\mathrm{HH}}\nolimits}

\begin{document}


\title{Singular Hochschild Cohomology and Gerstenhaber Algebra Structure}
\date{}

\author{Zhengfang WANG \thanks{
zhengfang.wang@imj-prg.fr, Universit\'e Paris Diderot-Paris 7, Institut
de Math\'ematiques de Jussieu-Paris Rive Gauche CNRS UMR 7586, B\^atiment Sophie Germain, Case 7012,
75205 Paris Cedex 13, France
}
}

\maketitle
\begin{abstract}
  In this paper, we define the singular Hochschild cohomology groups $\HH^i_{\sg}(A, A)$ of an associative $k$-algebra $A$ as morphisms from $A$ to $A[i]$ in the singular category $\DD_{\sg}(A\otimes_k A^{\op})$ for $i\in \Z$. We prove that $\HH^*_{\sg}(A, A)$ has a Gerstenhaber algebra structure and in the case of a symmetric algebra $A$, $\HH^*_{\sg}(A, A)$ is a Batalin-Vilkovisky (BV) algebra.
\end{abstract}

\section{Introduction}

Let $A$ be an associative algebra over a commutative ring $k$ such that $A$ is projective as a $k$-module. Then the Hochschild cohomology groups $\HH^i(A, A)$ can be defined as morphisms
from $A$ to $A[i]$ in the bounded derived category $\DD^b(A\otimes_k A^{\op})$ of the enveloping algebra $A\otimes_k A^{\op}$ for $i\in\Z_{\geq 0}$. Namely, we have
$$\HH^i(A, A):=\Hom_{\DD^b(A\otimes_k A^{\op})}(A, A[i]).$$
M. Gerstenhaber showed in \cite{Ger1} that there is a very rich structure on $\HH^*(A, A)$.
More precisely, he proved that $\HH^*(A, A)$ is a so-called Gerstenhaber algebra. Namely, there is a
Gerstenhaber bracket $[\cdot,\cdot]$ such that $[\cdot,\cdot]$ is a Lie bracket of degree $-1$, and a graded commutative associative cup product $\cup$, such that $[\cdot,\cdot]$ is a graded derivation of the cup product $\cup$ in each variable.

In this paper, we will generalize the Hochschild cohomology groups
to define the singular Hochschild cohomology groups $\HH_{\sg}^i(A, A)$ for
$i\in \Z$. Namely, we define the singular Hochschild cohomology groups as
$$\HH_{\sg}^i(A, A):=\Hom_{\DD_{\sg}(A\otimes_k A^{\op})}(A, A[i]),$$
where $\DD_{\sg}(A\otimes_k A^{\op})$ is the singular category of $A\otimes_k A^{\op}$.
Recall that $\DD_{\sg}(A\otimes_k A^{\op})$ is the Verdier quotient of $\DD^b(A\otimes_k A^{\op})$ by the full subcategory
consisting of perfect complexes, that is, bounded complexes of projective $A$-$A$-bimodules (cf. \cite{Bu, Orl}). We observe that $\HH^*_{\sg}(A, A)=0$ for an algebra $A$ of finite global dimension since the singular category
$\DD_{\sg}(A\otimes_k A^{\op})$ is zero in this case.
So, from this point of view, the algebras we are interested in are those
of infinite global dimension. Note that in general, $\HH^{i}_{\sg}(A, A)$
does not vanish even for $i\in \Z_{<0}.$
 The main result of this paper is as follow.
\begin{thm}[=Theorem \ref{thm-gerst}]
  Let $A$ be an associative algebra over a commutative ring $k$ such that $A$
  is projective as a $k$-module. Then the singular Hochschild cohomology
  $$\HH_{\sg}^*(A, A):=\bigoplus_{i\in\Z} \HH^i_{\sg}(A, A)$$
  is a  Gerstenhaber algebra, equipped with a Gerstenhaber bracket $[\cdot,\cdot]$ and
  the Yoneda product $\cup$ in the singular category $\DD_{\sg}(A\otimes A^{\op})$.
\end{thm}

From Buchweitz's work in his manuscript \cite{Bu}, we have a nice description on the singular Hochschild cohomology $\HH^*_{\sg}(A, A)$ for a self-injective algebra $A$. Namely, suppose that  $A$ is a self-injective algebra over a field $k$ (e.g. a group algebra of a finite group).
Then we have
\begin{equation*}
  \HH_{\sg}^i(A, A)\cong
  \begin{cases}
    \HH^i(A, A) &  \mbox{if}\ \  i>0,\\
    \HH_{-i-1}(A, \Hom_{A^e}(A, A^e)) & \mbox{if} \ i<-1.
  \end{cases}
\end{equation*}
Here we remark that, from \cite{CiSo}, for the case of a group algebra $k[G]$, where $G$ is a finite abelian group and the characteristic of $k$ divides the order of $G$,
$\HH_{\sg}^*(A, A)$ is very related to the Tate cohomology $\widehat{\HH}^*(G, k)$ of $G$ with coefficients in $k$, the trivial $kG$-module. We also remark that
in the case of a self-injective algebra $A$, the singular Hochschild cohomology $\HH_{\sg}^*(A, A)$ agrees with
the Tate-Hochschild cohomology defined in \cite{BeJo} and \cite{BeJoOp} and the stable Hochschild cohomology defined in \cite{EuSc}.

Recall that L. Menichi in \cite{Men} and T. Tradler in \cite{Tra} independently showed that the Hochschild cohomology $\HH^*(A, A)$ of a finite dimensional symmetric algebra $A$ has a new structure, the so-called Batalin-Vilkovisky (BV) structure (c.f.
Definition \ref{defn-BV}),
which has been studied in topology and mathematical physics during several
decades. Roughly speaking a BV structure is a differential operator
on Hochschild cohomology and it is a ``generator'' of the Gerstenhaber
bracket $[\cdot,\cdot]$, which means that $[\cdot,\cdot]$ is the obstruction of the
differential operator being a graded derivation with respect to the cup product. In this paper, we will generalize this result and  prove that
the singular Hochschild cohomology $\HH^*_{\sg}(A, A)$ of a finite dimensional symmetric algebra $A$ has a BV algebra structure. Namely, we have the following result.
\begin{thm}[=Corollary \ref{cor-bv}]
  Let $A$ be a symmetric algebra over a field $k$. Then the singular Hochschild cohomology $\HH^*_{\sg}(A, A)$ is a  BV algebra with the BV operator $\Delta_{\sg}$,  which is the Connes B-operator for the negative part $\HH^{< 0}_{\sg}(A, A)$, the $\Delta$-operator for the positive part $\HH^{> 0}_{\sg}(A, A)$ and $$\Delta_{\sg}|_{\HH^0_{\sg}(A, A)}=0: \HH^0_{\sg}(A, A)\rightarrow
  \HH^{-1}_{\sg}(A, A).$$ In particular, we have two BV subalgebras $\HH_{\sg}^{\leq 0}(A, A)$ and $\HH_{\sg}^{\geq 0}(A, A)$ with induced BV algebra structures.
\end{thm}
As a corollary, we obtain that the cyclic homology $\HC_*(A, A)[-1]$ of a symmetric algebra $A$, has a graded Lie algebra structure (cf. Corollary \ref{cor-cy}). We remark that
L. Menichi showed that the negative cyclic cohomology $\HC^*_{-}(A,A)[-1]$ (cf. Proposition 25, \cite{Men1}) and 
the cyclic cohomology $\HC_{\lambda}^*(A)[-1]$ (cf. Corollary 43, \cite{Men}) are both a graded Lie algebra . So in some sense our result is a dual
version of his results.


Throughout this paper, we fix a commutative ring $k$ with unit. We assume that all rings and the modules
are simultaneously $k$-modules and that all operations on rings and modules are naturally $k$-module homomorphisms. For simplicity,  we often use the symbol $\otimes$ to represent $\otimes_k,$ the tensor product over the commutative base ring $k$. For a $k$-algebra $A$, we denote $(a_i\otimes a_{i+1}\otimes\cdots\otimes a_{j})\in A^{\otimes j-i+1}(i\leq j)$ sometimes by $a_{i, j}$ for short, and denote the enveloping algebra $A\otimes_k A^{\op}$ by $A^e$.

\section*{Acknowledgement} This work is a part of author's PhD thesis, I would like to thank my PhD supervisor
Alexander Zimmermann for introducing this interesting topic and for his many valuable suggestions for improvement. I am grateful to Claude Cibils and
Selene Sanchez for some interesting discussions,  to Murray Gerstenhaber for some remarks on his paper \cite{Ger1}  and to Reiner Hermann for
some discussions on his PhD thesis when I just started this project.
I also would like to thank Ragnar-Olaf Buchweitz for useful suggestions during this project. Special thanks to my PhD co-supervisor, Marc Rosso for
his constant support and encouragement during my career in mathematics.

\section{Preliminaries}
In this section we recall some notions on Hochschild cohomology and  Gerstenhaber algebras. For more details, we refer the reader to \cite{Lod, Ger1}.
\begin{defn}
Let $k$ be a commutative ring with unit. A differential graded Lie algebra (DGLA) is a differential $\Z$-graded $k$-module $(L, d)$ with a bracket $[\cdot,\cdot]: L^i\times
L^j\rightarrow L^{i+j}$ which satisfies the following properties:
\begin{enumerate}
\item it is skew-symmetric: $$[\alpha, \beta]=-(-1)^{|\alpha||\beta|} [\beta, \alpha];$$
\item satisfies the graded Leibniz rule: $$d([\alpha, \beta])=(-1)^{|\beta|}[d\alpha, \beta]+[\alpha, d\beta];$$
\item and the graded Jacobi identity: $$(-1)^{(|\alpha|-1)(|\gamma|-1)}[[\alpha, \beta], \gamma]+(-1)^{(|\beta|-1)(|\alpha|-1)}[[\beta,\gamma],\alpha]+(-1)^{(|\gamma|-1)(|\beta|-1)}[[\gamma, \alpha], \beta]=0,$$
\end{enumerate}
where $\alpha, \beta, \gamma$ are arbitrary homogeneous elements in $(L, d)$ and $|\alpha|$ is the degree of the homogeneous element $\alpha$.
\end{defn}
\begin{rem}
Let $(L, d, [\cdot,\cdot])$ be a DGLA. Then the homology $H^*(L, d)$ of the differential graded module $(L, d)$ is a $\Z$-graded Lie algebra with the induced bracket $[\cdot,\cdot]$.
\end{rem}
\begin{defn}
A Gerstenhaber algebra is a $\Z$-graded $k$-module $\calH^*:=\oplus_{n\in\Z}\calH^n$ equipped with:
\begin{enumerate}
  \item a graded commutative associative product $\cup$ of degree zero, with unit $1\in\calH^0$,
  \begin{eqnarray*}
    \begin{tabular}{cccc}
    $\cup:$ & $\calH^m \times \calH^n$ & $\rightarrow$ & $\calH^{m+n}$\\
    & $(\alpha, \beta)$ & $\mapsto$ & $\alpha\cup\beta.$
    \end{tabular}
  \end{eqnarray*}
  In particular, $\alpha\cup \beta=(-1)^{|\alpha||\beta|}\beta\cup\alpha$;
  \item a graded Lie algebra structure $[\cdot, \cdot]$ on $\calH^*[-1]$, that is,
  $$[\alpha, \beta]=-(-1)^{(|\alpha|-1)(|\beta|-1)}[\beta, \alpha]$$
  and
  $$(-1)^{(|\alpha|-1)(|\gamma|-1)}[[\alpha, \beta], \gamma]+(-1)^{(|\beta|-1)(|\alpha|-1)}[[\beta,\gamma],\alpha]+(-1)^{(|\gamma|-1)(|\beta|-1)}[[\gamma, \alpha], \beta]=0;$$
  \item compatibility between $\cup$ and $[\cdot, \cdot]$:
  $$[\alpha, \beta\cup \gamma]=\beta\cup[\alpha, \gamma]+(-1)^{|\gamma|(|\alpha|-1)}[\alpha, \beta]\cup \gamma,$$
  (or equivalently,
  $$[\alpha\cup \beta, \gamma]=[\alpha, \gamma]\cup \beta+(-1)^{|\alpha|(|\gamma|-1)}\alpha\cup[\beta, \gamma])$$
\end{enumerate}
where $\alpha, \beta, \gamma$ are arbitrary homogeneous elements in $\calH^*$ and $|\alpha|$ is the degree of the homogeneous element $\alpha$.
\end{defn}


We follow \cite{BaGi} to define a Gerstenhaber module of a Gerstenhaber algebra $\calH^*$.
\begin{defn}\label{defn-module}
  A Gerstenhaber module of $\calH^*$ is a $\Z$-graded vector space $\calF^*$
  equipped with:
  \begin{enumerate}
    \item a module structure $\cup'$ for the graded algebra $(\calH^*, \cup)$;
    \item a module structure $[\cdot, \cdot]'$ for the graded Lie algebra $(\calH^*[-1], [\cdot, \cdot])$. That is,
    $$[[\alpha, \beta], x]'=[\alpha, [\beta, x]']'-(-1)^{(|\alpha|-1)(|\beta|-1)}[\beta,[\alpha, x]']';$$
    \item compatibility:
      $$[\alpha\cup\beta, x]'=(-1)^{|\alpha|(|x|-1)}\alpha\cup'[\beta, x]'+(-1)^{(|\alpha|+|x|-1)|\beta|} \beta\cup'[\alpha, x]',$$
      $$[\alpha, \beta\cup'x]'=\beta\cup'[\alpha, x]'+(-1)^{|x|(|\alpha|-1)}[\alpha, \beta]\cup'x
    $$
    where $\alpha, \beta$ are arbitrary homogeneous elements in $\calH^*$ and $x$ is arbitrary homogeneous element in $\calF^*$.
  \end{enumerate}
\end{defn}
\begin{rem}
For any Gerstenhaber algebra $\calH^*$, $\calH^*$ is a Gerstenhaber module over itself.
\end{rem}
Naturally, we can define morphisms between Gerstenhaber modules.
\begin{defn}
Let $(\calH^*, \cup, [\cdot, \cdot])$ be a Gerstenhaber algebra over a commutative ring $k$. Let $(\calF_1^*, \cup_1,[\cdot, \cdot]_1)$ and $(\calF^*_2, \cup_2, [\cdot,\cdot]_2)$ be two Gerstenhaber modules of $\calH^*$.
We say that a $k$-module morphism $\varphi: \calF^*_1\rightarrow \calF^*_2$ is a Gerstenhaber morphism of degree $r, (r\in\Z)$ if the following two conditions are satisfied
\begin{enumerate}
\item $\varphi$ is a module morphism of degree $r$ for the graded commutative algebra  $(\calH^*, \cup)$. That is, for any
$f\in \calH^m$ and $g\in \calF_1^*$, $$\varphi(f\cup_1 g)=(-1)^{mr}f\cup_2\varphi(g)$$
\item $\varphi$ is a module morphism of degree $r$ for the graded Lie algebra  $(\calH^*[-1], [\cdot, \cdot]).$ That is,
for any $f\in \calH^m$ and $g\in \calF_1^*$,
$$\varphi([f, g]_1)=(-1)^{(m-1)r}[f, \varphi(g)]_2.$$
\end{enumerate}
\end{defn}

The classical example of Gerstenhaber algebras is the Hochschild cohomology $$\HH^*(A, A):=\bigoplus_{n\in\Z_{\geq 0}}\HH^n(A, A)$$ of an associative algebra $A$. Let us start to recall some notions on Hochschild cohomology.
Let $A$ be an associative algebra over $k$ such that $A$ is projective as a $k$-module and $M$ be an $A$-$A$-bimodule. Recall that the Hochschild cohomology of $A$ with coefficients in $M$ is defined as
$$\HH^*(A, M):=\Ext_{A^e}^*(A, M),$$
and Hochschild homology is defined as
$$\HH_*(A, M):=\Tor_*^{A^e}(A, M)$$
where $A^e:=A\otimes_k A^{\op}$ is the enveloping algebra of $A$. Recall that we have the following (un-normalized) bar resolution of $A$,
\begin{equation}\label{bar}
\xymatrix{
\Barr_*(A): \cdots\ar[r] &  A^{\otimes(r+2)}\ar[r]^{d_r} & A^{\otimes(r+1)} \ar[r] & \cdots \ar[r] & A^{\otimes 3} \ar[r]^{d_1} & A^{\otimes 2} \ar[r]^{d_0:=\mu} & A
},
\end{equation}
where $\mu$ is the multiplication of $A$ and $d_r$ is defined as follows,
$$d_r(a_0\otimes a_1\otimes \cdots \otimes a_{r+2})=\sum_{i=0}^r(-1)^ia_{0, i-1}\otimes a_ia_{i+1}\otimes a_{i+2, r+1}.$$
Denote
\begin{equation*}
\begin{split}
  C^r(A, M):&=\Hom_{A^e}(A^{\otimes r+2}, M), \\
  C_r(A, M):&= M\otimes_{A^e}A^{\otimes r+2},
\end{split}
\end{equation*}
for any $r\in\Z_{\geq 0}$. Note that
\begin{equation*}
  \begin{split}
    \Hom_{A^e}(A^{\otimes r+2}, M)&\cong \Hom_{k}(A^{\otimes r}, M),\\
     M\otimes_{A^e}A^{\otimes r+2}&\cong  M\otimes_k A^{\otimes r}.
  \end{split}
\end{equation*}
We also consider the the normalized bar resolution $\overline{\Barr}_*(A)$, which is defined as $$\overline{\Barr}_r(A):=A\otimes \overline{A}^{\otimes r}\otimes A,$$
where $\overline{A}:=A/(k \cdot 1_A)$, with induced differential in $\Barr_*(A)$.
Thus the Hochschild cohomology $\HH^*(A, M)$ can be computed by the following complex,
$$\xymatrix{
C^*(A, M): M\ar[r]^-{\delta^0} &  \Hom_k(A, M)\ar[r] & \cdots \ar[r] & \Hom_k(A^{\otimes r}, M) \ar[r]^{\delta^r} & \Hom_k(A^{\otimes r+1}, M)  \ar[r] &\cdots},$$
where $\delta^r$ is defined as follows, for any $f\in \Hom_k(A^{\otimes r}, M)$
\begin{eqnarray*}
  \delta^r(f)(a_1\otimes\cdots \otimes a_{r+1}):&=&a_1f(a_{2,r+1})+\sum_{i=1}^r(-1)^if(a_{1,i-1}\otimes a_ia_{i+1}\otimes a_{i+2, r+1})+\\
  &&(-1)^{r+1}f(a_{1,r})a_{r+1}.
\end{eqnarray*}
We denote $$Z^r(A, M):=\Ker(\delta^r)$$ and $$B^r(A, M):=\Imm(\delta^{r-1})$$ for any $r\in \Z_{\geq 0}$.
Then we have $$\HH^r(A, M)=Z^r(A, M)/B^r(A, M).$$
We can also compute the Hochschild cohomology $\HH^*(A,A)$ via the normalized Bar resolution $\overline{\Barr}_*(A)$, namely, we have (cf. e.g. \cite{Lod})
$$\HH^r(A, M)=\overline{Z}^r(A, M)/\overline{B}^r(A, M),$$
where $\overline{Z}^r(A, M)$ and $\overline{B}^r(A, M)$ are respectively the r-th cocycle and r-th coboundary in the normalized cochain complex $\overline{C}^*(A, M)$.

The Hochschild homology $\HH_*(A, M)$ is the homology of the following complex,
$$\xymatrix{
C_*(A, M): \cdots\ar[r] &  C_r(A, M) \ar[r]^{\partial_r} & C_{r-1}(A, M)\ar[r]  & \cdots \ar[r] & C_0(A, M)
}
$$
where
\begin{equation*}
  \begin{split}
    \partial_r(m\otimes a_1\otimes\cdots\otimes a_r):=
    ma_1\otimes a_{2,r}+\sum_{i=1}^{r-1}(-1)^{i}m\otimes a_{1, i-1}\otimes a_ia_{i+1}\otimes a_{i+2, r}+(-1)^ra_rm\otimes a_{1, r-1}.
  \end{split}
\end{equation*}
Denote
\begin{equation*}
  \begin{split}
    B_r(A, M):&=\Imm(\partial_{r+1})\\
    Z_r(A, M):&=\Ker(\partial_r).
  \end{split}
\end{equation*}
Then we have
\begin{equation*}
  \HH_r(A, M):=Z_r(A, M)/B_r(A, M).
\end{equation*}
Similarly, it also can be computed by the normalized bar resolution $\overline{\Barr}_*(A, M)$, namely, we have (cf. e.g. \cite{Lod})
$$\HH_r(A, M)=\overline{Z}_r(A, M)/\overline{B}_r(A, M).$$
Let us recall the cup product,
$$\cup: C^m(A, A)\times C^n(A, M)\rightarrow C^{m+n}(A, M),$$
which is defined in the following way. Given $f\in C^m(A, A)$ and $g\in C^n(A, M)$,
$$(f\cup g)(a_1\otimes \cdots \otimes a_{m+n}):=f(a_{1, m})g(a_{m+1,m+n}).$$
One can check that this cup product $\cup$ induces a well-defined operation (still denoted by $\cup$)
on cohomology groups, that is,
$$\cup: \HH^m(A, A)\times \HH^n(A, M)\rightarrow \HH^{m+n}(A, M).$$
Recall that there is also a circ product,
$$\circ: C^n(A, M)\times C^m(A, A)\rightarrow C^{m+n-1}(A, M)$$
which is defined as follows, given $f\in C^m(A, A)$ and $g\in C^n(A, M)$,
for $1\leq i\leq n$, set
$$g\circ_if(a_1\otimes \cdots\otimes a_{m+n-1}):=g(a_{1, i-1} \otimes f(a_{i, i+m-1})\otimes a_{i+m,m+n-1}),$$
\begin{equation}\label{equ-circ-product}
  g\circ f:=\sum_{i=1}^n(-1)^{(m-1)(i-1)}g\circ_if,
\end{equation}
for $n=0$, we set $$g\circ f:=0.$$
Using this circ product, one can define a Lie bracket  $[\cdot, \cdot]$ on $\HH^*(A, A)$ in the following way. Let $f\in C^m(A, A)$ and $g\in C^n(A, A)$,
define $$[f, g]:=f\circ g-(-1)^{(m-1)(n-1)}g\circ f.$$
One can check that this Lie bracket induces a well-defined Lie bracket (still denoted by $[\cdot, \cdot]$) on $\HH^*(A, A)$.

With these two operators $\cup$ and $[\cdot, \cdot]$ on $\HH^*(A, A)$, Gerstenhaber proves the following result.
\begin{thm}[\cite{Ger1}]
Let $A$ be an associative algebra over a commutative ring $k$. Then
 the Hochschild cohomology $\HH^*(A, A)$, equipped with the cup product $\cup$ and bracket $[\cdot,\cdot]$ is a Gerstenhaber algebra.
\end{thm}
\begin{rem}
  Let $A$ be an associative $k$-algebra such that $A$ is projective as a $k$-module.
  Then for any $m\in \Z_{\geq 0}$, $$\HH^m(A, A)\cong\Hom_{\DD^b(A^e)}(A, A[m]),$$ and the cup product $\cup$ can
  be interpreted as compositions of morphisms (namely, the Yoneda product) in $\DD^b(A\otimes_k A^{\op})$.
\end{rem}

At the end of this section, let us recall the cap product $\cap$, which is an action of Hochschild cohomology on Hochschild homology.
For any $r, p\in\Z_{\geq 0}$ such that $r\geq p$, there is a bilinear map
\begin{equation*}
  \cap: C_r(A, M)\otimes C^p(A, A)\rightarrow C_{r-p}(A, M)
\end{equation*}
sending $(m\otimes a_1\otimes \cdots \otimes a_r)\otimes \alpha$
to
\begin{equation*}
  (m\otimes a_1\otimes \cdots \otimes a_r)\cap \alpha:=
  (-1)^{rp}(m\otimes_A\alpha(a_1\otimes \cdots \otimes a_p)\otimes
  a_{p+1}\otimes \cdots \otimes a_{r})
\end{equation*}
It is straightforward to verify that $\cap$ induces a well-defined map, which we still denote by $\cap$, on the
level of homology,
\begin{equation}\label{equ-cap}
  \cap: \HH_r(A, M)\otimes \HH^p(A, A)\rightarrow \HH_{r-p}(A, M).
\end{equation}

\section{Singular Hochschild cohomology}

  Let $A$ be an associative algebra over a commutative ring $k$ such that $A$ is projective as a $k$-module.
 Recall the un-normalized bar resolution $\Barr_*(A)$ (cf. (\ref{bar})) of $A$,
\begin{eqnarray*}
  \xymatrix{
 \cdots\ar[r]^{d_2}& A^{\otimes 3} \ar[r]^-{d_1} & A^{\otimes 2} \ar[r]^{d_0} & A\ar[r] &0.
  }
\end{eqnarray*}
Let us denote the $p$-th kernel $\Ker(d_{p-1})$ in the un-normalized bar resolution $\Barr_*(A)$ by $\Omega^p(A)$.
Then we have the following short exact sequence for $p\in\Z_{>0}$,
\begin{eqnarray}\label{short-exact}
  0\rightarrow\Omega^p(A)\rightarrow A^{\otimes (p+1)} \rightarrow \Omega^{p-1}(A)\rightarrow 0
\end{eqnarray}
which induces a long exact sequence
\begin{eqnarray}\label{longexact}
  \cdots\rightarrow  \HH^m(A, A^{\otimes (p+1)}) \rightarrow \HH^m(A, \Omega^{p-1}(A))\rightarrow \HH^{m+1}(A, \Omega^p(A))\rightarrow\cdots
\end{eqnarray}
We denote the connecting morphism in (\ref{longexact}) by,
$$\theta_{m,p-1}: \HH^m(A, \Omega^{p-1}(A))\rightarrow \HH^{m+1}(A, \Omega^p(A))$$
for $m\in\Z_{\geq 0}$.
Hence we obtain an inductive system for any fixed $m\in\Z_{\geq 0}$,
\begin{eqnarray}\label{system}
  \xymatrix{
  \HH^m(A, A)\ar[r]^-{\theta_{m, 0}} & \HH^{m+1}(A, \Omega^1(A)) \ar[r]^{\theta_{m,1}} & \HH^{m+2}(A, \Omega^2(A))\ar[r] &\cdots
  }
\end{eqnarray}
Let
$$\lim_{\substack{\longrightarrow\\p\in\Z_{\geq 0}}} \HH^{m+p}(A, \Omega^p(A))$$
be the colimit of the inductive system (\ref{system}) above.

Since $A$ is a $k$-algebra such that $A$ is projective as a $k$-module,
we have a canonical isomorphism $$\HH^{m}(A, M)\cong \Hom_{\DD^b(A^e)}(A, M[m])$$
for any $A$-$A$-bimodule $M$ and $m\in \Z_{\geq 0}$.
Hence we have the following morphism for any $m\in \Z, p\in \Z_{\geq 0}$ such that $m+p>0$,
\begin{equation}\label{equ-natural}
  \Phi_{m, p}:\HH^{m+p}(A, \Omega^p(A))\rightarrow \Hom_{\DD_{\sg}(A^e)}(A, \Omega^p(A)[m+p])\rightarrow \Hom_{\DD_{\sg}(A^e)}(A, A[m]),
\end{equation}
which are compatible with the inductive system (\ref{system}) above.
So the collection of maps $\Phi_{m, p}$ induces a morphism for any fixed $m\in \Z$,
$$\Phi_m: \lim_{\substack{\longrightarrow\\p\in\Z_{\geq 0}\\ m+p>0}} \HH^{m+p}(A, \Omega^p(A))\rightarrow \Hom_{\DD_{\sg}(A^e)}(A, A[m]).$$

Next we will prove that $\Phi_m$ is an isomorphism for any $m\in\Z$.
\begin{prop}\label{prop}
  For any $m\in\Z$, the morphism $$\Phi_m: \lim_{\substack{\longrightarrow\\p\in\Z_{\geq 0}\\ m+p>0}} \HH^{m+p}(A, \Omega^p(A))\rightarrow \Hom_{\DD_{\sg}(A^e)}(A, A[m])$$ defined above is an isomorphism.
\end{prop}
\pf First, let us recall the following fact (cf. e.g. Proposition 6.7.17 \cite{Zim}):
\begin{fact}
  The following canonical homomorphism is an isomorphism for any $m\in\Z$
$$\Hom_{\DD_{\sg}(A^e)}(A, A[m])\cong \lim_{\substack{\longrightarrow\\p\in\Z_{\geq 0}\\m+p>0}} \underline{\Hom}_{A^e}(\Omega^{m+p}(A), \Omega^p(A)).$$
\end{fact}
Now using the fact above, we obtain that $\Phi_m$ is surjective. Indeed,
assume $$f\in \Hom_{\DD_{\sg}(A^e)}(A, A[m]),$$ then from the fact above,  there exists $p\in\Z_{\geq 0}$ such that $f$ can be represented by some element
$$f'\in \underline{\Hom}_{A^e}(\Omega^{m+p}(A), \Omega^p(A)),$$
hence $f$ is also represented by some element $f''\in \Hom_{A^e}(\Omega^{m+p}(A), \Omega^p(A)).$
It follows that $f''$ induces a cocycle (See Diagram (\ref{diagram2})) $$\alpha:=f''\circ d_{m+p}\in \Hom_{A^e}(A^{\otimes m+p+2}, \Omega^p(A)),$$
hence $\alpha\in \HH^{m+p}(A, \Omega^p(A))$.
\begin{eqnarray}\label{diagram2}
  \xymatrix{
  A^{\otimes m+p+3} \ar[d]_{d_{m+p+1}}\\
  A^{\otimes m+p+2} \ar[r]^{\alpha}\ar@{->>}[d]_{d_{m+p}}  & \Omega^p(A)\\
  \Omega^{m+p}(A)\ar[ru]_{f''}\ar@{_(->}[d]\\
  A^{\otimes m+p+1}
  }
\end{eqnarray}
Moreover, we have that $\Phi_m(\alpha)=f$, so $\Phi_m$ is surjective.
Here we remark that
the following morphism between two inductive systems,
\begin{eqnarray*}
  \xymatrix{
   \cdots\ar[r]&\HH^{m+p}(A, \Omega^p(A))\ar[r]\ar[d] & \HH^{m+p+1}(A,\Omega^{p+1}(A))\ar[r]\ar[d]& \cdots\\
    \cdots\ar[r] & \underline{\Hom}_{A^e}(\Omega^{m+p}(A), \Omega^p(A))\ar[r] & \underline{\Hom}_{A^e}(\Omega^{m+p+1}(A), \Omega^{p+1}(A))\ar[r]&\cdots
  }
\end{eqnarray*}
induces a commutative diagram,
$$
  \xymatrix{
  \lim\limits_{\substack{\longrightarrow\\p\in\Z_{\geq 0}\\ m+p>0}} \HH^{m+p}(A, \Omega^p(A))\ar[r]^-{\Phi_m}\ar[d] &  \Hom_{\DD_{\sg}(A^e)}(A, A[m])\\
  \lim\limits_{\substack{\longrightarrow\\p\in\Z_{\geq 0}\\m+p>0}} \underline{\Hom}_{A^e}(\Omega^{m+p}(A), \Omega^p(A))\ar[ru]_{\cong}.
  }
$$

It remains to prove that $\Phi_m$ is injective.
Assume that there exists $$\beta\in \lim_{\substack{\longrightarrow\\p\in\Z_{\geq 0}\\ m+p>0}} \HH^{m+p}(A, \Omega^p(A))$$
such that $\Phi_m(\beta)=0$. Then there exists $p\in\Z_{\geq 0}$ such that
$\beta$ is represented by some element $\beta'\in \HH^{m+p}(A, \Omega^p(A))$ and $\beta'$ is mapped into zero
under the morphism $$\HH^{m+p}(A, \Omega^p(A))\rightarrow \underline{\Hom}_{A^e}(\Omega^{m+p}(A), \Omega^p(A)).$$
So from this it follows that the cocycle $\beta'$ induces a morphism $g':\Omega^{m+p}(A)\rightarrow \Omega^p(A)$ such that $g'$ factors through a projective $A$-$A$-bimodule $P$.
The maps are illustrated by the following diagram.
\begin{equation}\label{diagram-tau}
  \xymatrix{
  A^{\otimes m+p+3} \ar[d]_{d_{m+p+1}}\\
  A^{\otimes m+p+2} \ar[r]^{\beta'}\ar@{->>}[d]_{d_{m+p}}  & \Omega^p(A)\\
  \Omega^{m+p}(A)\ar[ru]_-{g'}\ar@{_(->}[d] \ar[r]_-{\sigma}  & P\ar[u]_-{\tau}\\
  A^{\otimes m+p+1}
  }
\end{equation}
By funtctoriality, the $A$-$A$-bimodule morphism $\tau:P\rightarrow \Omega^p(A)$ in Diagram (\ref{diagram-tau}) induces the following map.
\begin{eqnarray}\label{tau}
  \begin{tabular}{ccccc}
$\tau^*$ : & $\HH^{m+p}(A, P)$&$\rightarrow$ & $\HH^{m+p}(A, \Omega^p(A))$\\
 &$\sigma$&$\mapsto$&$ \beta'$\\
\end{tabular}
\end{eqnarray}
Since $P$ is a projective $A$-$A$-bimodule, we have the following commutative diagram between two short exact sequences.
\begin{eqnarray*}
\xymatrix{
  0\ar[r] &  \Omega^{p+1}\ar[r] &   A^{p+2} \ar[r] &   \Omega^p(A)\ar[r] & 0\\
  0\ar[r] & 0\ar[r]  \ar[u]&   P\ar[r]^{\cong}\ar[u]^{\tau'} & P \ar[r]\ar[u]^-{\tau} & 0
  }
\end{eqnarray*}
Hence from the functoriality of long exact sequences induced from short exact sequences, we have the following commutative diagram.
\begin{eqnarray}\label{diagram1}
\xymatrix@C=1.8em{
   \cdots\ar[r]  &  \HH^{m+p}(A, A^{\otimes p+2}) \ar[r]  &  \HH^{m+p}(A, \Omega^{p}(A))\ar[r]^-{\theta_{m+p,p}} & \HH^{m+p+1}(A, \Omega^{p+1}(A))\\
   \cdots\ar[r] & \HH^{m+p}(A, P) \ar[u]^{\tau'^*} \ar[r]^{\cong} & \HH^{m+p}(A, P)\ar[u]^{\tau^*}\ar[r] & 0\ar[u]
   }
\end{eqnarray}
So from Diagram (\ref{diagram1}) above, we obtain that $$\theta_{m+p,p}(\beta')=\theta_{m+p, p}\tau^*(\sigma)=0,$$
thus it follows that $\beta$ is represented by zero element in $\HH^{m+p+1}(A, \Omega^{p+1}(A))$, hence $$\beta=0\in \lim_{\substack{\longrightarrow\\p\in\Z_{\geq 0}\\m+p>0}} \HH^{m+p}(A, \Omega^p(A)).$$ So  $\Phi_m$ is injective.
Therefore, $\Phi_m$ is an isomorphism.
\epf
\begin{rem}
 By the same argument,  we also have the following isomorphism for any $m\in \Z$,
  $$\overline{\Phi}_m: \lim_{\substack{\longrightarrow\\p\in\Z_{\geq 0}\\ m+p>0}} \HH^{m+p}(A, \overline{\Omega}^p(A))\rightarrow \Hom_{\DD_{\sg}(A^e)}(A, A[m])$$
  where $\overline{\Omega}^p(A)$ is the $p$-th kernel in the normalized bar resolution of $A$.
\end{rem}

\section{Gerstenhaber algebra structure on singular Hochschild cohomology}
In this section, we will prove the following main theorem.
\begin{thm}\label{thm-gerst}
Let $A$ be an associative algebra over a commutative ring $k$ and suppose that $A$ is projective as a $k$-module.
Then the singular Hochschild cohomology (shifted by $[1]$) $$\HH^*_{\sg}(A, A)[1]:=\bigoplus_{n\in\Z}\Hom_{\DD_{\sg}(A^e)}(A, A[n])[1]$$ is a Gerstenhaber algebra.
\end{thm}
For $p\in \Z_{\geq 0}$, we denote the $p$-th kernel $\Ker(\overline{d}_{p-1})$ in the normalized bar resolution
$\overline{\Barr}_*(A)$  by
$\overline{\Omega}^p(A).$
Note that $\overline{\Omega}^0(A)=A$.

Let $m, n\in \Z_{>0}$ and $p, q\in \Z_{\geq0}$.
we shall define a Gerstenhaber bracket as follows
\begin{equation}\label{equ-bracket}
  [\cdot,\cdot]: C^m(A, \overline{\Omega}^{p}(A))\otimes C^n(A, \overline{\Omega}^q(A))\rightarrow C^{m+n-1}(A, \overline{\Omega}^{p+q}(A)).
\end{equation}
Let $f\in C^m(A, \overline{\Omega}^p(A))$ and $g\in C^n(A, \overline{\Omega}^q(A))$,
define
\begin{equation*}
f\bullet_i g:=
\begin{cases}
  d((f\otimes \id^{\otimes q})(\id^{\otimes i-1}\otimes g\otimes \id^{\otimes m-i})\otimes 1)&\mbox{if} \ 1\leq i\leq m, \\
 d((\id^{\otimes -i}\otimes f\otimes \id^{\otimes q+i})(g\otimes \id^{\otimes m-1})\otimes 1) & \mbox{if} \ -q\leq i \leq -1,
\end{cases}
\end{equation*}
and
\begin{equation*}
  f\bullet g:=\sum_{i=1}^m(-1)^{r(m,p;n,q, i)} f\bullet_i g+\sum_{i=1}^q (-1)^{s(m,p; n, q; i)}f\bullet_{-i} g
\end{equation*}
where $r(m, p;n, q;i)$ and $s(m,p;n,q;i)$ are defined in the following way,
\begin{equation}\label{equ-coe}
\begin{split}
r(m, p;n, q;i)&=p+q+(i-1)(q-n-1), \ \ 1\leq i\leq m,\\
  s(m, p;n, q;i)&=p+q+i(p-m-1), \ \ 1\leq i\leq q.
\end{split}
\end{equation}
Then we define
\begin{equation}\label{equ-defn-bracket}
  [f, g]:=f\bullet g-(-1)^{(m-p-1)(n-q-1)} g\bullet f.
\end{equation}
\begin{rem}
 We have a double complex $C^{*, *}(A, A)$, which is defined by
 \begin{equation*}
 C^{m,p}(A, A):=
 \begin{cases}
 C^{m}(A, \overline{\Omega}^p(A)) & \mbox{if}\ \  m\in\Z_{>0}, p\in\Z_{\geq 0},\\
 0 & \mbox{otherwise}
 \end{cases}
 \end{equation*}
 with the horizontal differential $\delta: C^{m}(A, \overline{\Omega}^p(A))\rightarrow C^{m+1}(A, \overline{\Omega}^p(A))$, induced from bar resolution and the vertical differential zero. Recall that the total complex of $C^{*, *}(A, A)$, denoted by
 $\Tot (C^{*, *}(A, A))^*$, is defined as follows,
 \begin{equation*}
   \Tot(C^{*, *}(A, A))^n:=\bigoplus_{n=m-p} C^{m,p}(A,A)
 \end{equation*}
 with the differential induced from $C^{*, *}(A, A)$.
\end{rem}

Note that the bullet product $\bullet$ is not associative, however it has the following ``weak associativity''.
\begin{lemma}\label{lemma-bullet1}
Let $f_i\in C^{m_i}(A, \overline{\Omega}^{p_i}(A))$ for $i=1, 2, 3.$
\begin{enumerate}
  \item For $1\leq j\leq m_1$ and $1\leq i\leq m_1+m_2-1$, we have
  \begin{equation*}
    (f_1\bullet_j f_2)\bullet_i f_3=
    \begin{cases}
  (-1)^{p_1+p_3}f_1\bullet_j(f_2\bullet_{i-j+1}f_3)  &\mbox{if} \ \  0\leq i-j< m_2, \\
  (-1)^{p_3+p_2}(f_1\bullet_{i+p_2-m_2+1} f_3)\bullet_j f_2 & \mbox{if} \ \ m_2\leq i-j, i<m_1+m_2-p_2\\
  f_3\bullet_{-(p_1+p_2+1+i-m_1-m_2)}(f_1\bullet_j f_2) &\mbox{if} \ \
  m_2\leq i-j, m_1+m_2-p_2\leq i\\
  (-1)^{p_1+p_3}f_1\bullet_i (f_2\bullet_{-(j-i)} f_3)  & \mbox{if} \ \ 1\leq j-i\leq p_3 \\
  (-1)^{p_2+p_3}(f_1\bullet_if_3)\bullet_{m_3+j-p_3-1} f_2 & \mbox{if}
  \ \ p_3<j-i\\
\end{cases}
  \end{equation*}
  \item For $1\leq j\leq p_2$ and $1\leq i\leq m_1+m_2-1$,
  \begin{equation*}
    (f_1\bullet_{-j} f_2)\bullet_if_3=
    \begin{cases}
      (-1)^{p_1+p_3} f_1\bullet_{-j}(f_2\bullet_if_3) & \mbox{if} \ \ 1\leq i\leq m_2\\
      (-1)^{p_2+p_3}(f_1\bullet_{p_2-j+i-m_2+1} f_3)\bullet_{-j} f_2 & \mbox{if} \ \ m_2<i<m_1+m_2-p_2+j\\
      f_3\bullet_{-(p_1+p_2+1+i-m_1-m_2)}(f_1\bullet_{-j} f_2) &
      \mbox{if} \ \ m_2<m_1+m_2-p_2+j\leq i\\
       f_3\bullet_{-(p_1+p_2+1+i-m_1-m_2)}(f_1\bullet_{-j} f_2) &
      \mbox{if} \ \ m_2<i, m_1+j\leq p_2\\
    \end{cases}
  \end{equation*}
  \item For $1\leq j\leq m_1$ and $1\leq i\leq p_3$,
  \begin{equation*}
  (f_1\bullet_jf_2)\bullet_{-i} f_3=
    \begin{cases}
   (-1)^{p_1+p_3} f_1\bullet_{-i}(f_2\bullet_{-(i+j-1)} f_3) & \mbox{if}\ \ i+j<p_3-2\\
      (-1)^{p_2+p_3}(f_1\bullet_{-i}f_3)\bullet_{m_3-p_3+j+i-1} f_2 & \mbox{if} \ \ p_3-2\leq i+j
    \end{cases}
  \end{equation*}
  \item For $1\leq j\leq p_2$ and $1\leq i\leq p_3$,
  \begin{equation*}
    (f_1\bullet_{-j} f_2)\bullet_{-i} f_3=(-1)^{p_1+p_3}f_1\bullet_{-(i+j)}(f_2\bullet_{-i} f_3).
  \end{equation*}
  \end{enumerate}
  \end{lemma}
  Similarly, we have the following lemma.
 \begin{lemma}\label{lemma-bullet2}
 \begin{enumerate}
 \item For $1\leq j\leq m_1$ and $1\leq i\leq m_2$,
 \begin{equation*}
 f_1\bullet_j(f_2\bullet_i f_3)=(f_1\bullet_jf_2)\bullet_{j+i-1} f_3
 \end{equation*}
  \item For $1\leq j\leq m_1$ and $1\leq i\leq p_3$,
  \begin{equation*}
  f_1\bullet_j(f_2\bullet_{-i} f_3)=
   \begin{cases}
   (-1)^{p_1+p_3}(f_1\bullet_{i+j} f_2)\bullet_{j} f_3 & \mbox{if} \ \ i+j\leq m_1\\
  (-1)^{p_1+p_2}f_2\bullet_{-(i+j-m_1+p_1)}(f_1\bullet_j f_3) & \mbox{if}\ \
  m_1<i+j
      \end{cases}
  \end{equation*}
  \item For $1\leq j\leq p_2+p_3$ and $1\leq i\leq m_2$,
  \begin{equation*}
  f_1\bullet_{-j}(f_2\bullet_if_3)=
  \begin{cases}
  (-1)^{p_1+p_3}(f_1\bullet_{-j} f_2)\bullet_i f_3 & \mbox{if} \ \  0\leq j\leq p_2\\
  (-1)^{p_1+p_2}f_2\bullet_i(f_1\bullet_{-(j+m_2-p_2-i)}f_3) & \mbox{if}\ \ 0<j-p_2\leq p_3+i-m_2\\
 (f_2\bullet_if_3)\bullet_{m_2+m_3-p_2-p_3-1+j} f_1  & \mbox{if}\ \ 1\leq p_3+i-m_2+1\leq j-p_2\\
 (f_2\bullet_if_3)\bullet_{m_2+m_3-p_2-p_3-1+j} f_1 &  \mbox{if}\ \
 1\leq j-p_2,  p_3+i-m_2<0
  \end{cases}
  \end{equation*}
  \item For $1\leq j\leq p_2+p_3$ and $1\leq i\leq p_3$,
  \begin{equation*}
  f_1\bullet_{-j}(f_2\bullet_{-i} f_3)=
  \begin{cases}
 (-1)^{p_1+p_3} (f_1\bullet_{i-j+1} f_2)\bullet_{-j} f_3 & \mbox{if} \ \ 0\leq i-j\\
  (-1)^{p_1+p_3}(f_1\bullet_{-(j-i)} f_2)\bullet_{-i} f_3 & \mbox{if} \ \ 0<j-i\leq p_2\\
  (f_2\bullet_{-i} f_3)\bullet_{m_2+m_3-p_2-p_3-1+j}f_1 &\mbox{if} \ \ p_2<j-i
  \end{cases}
  \end{equation*}
    \end{enumerate}
\end{lemma}

\begin{rem}\label{rem-delta}
  Similar to \cite{Ger1}, the cup product $\cup$ for $C^{*,*}(A, A)$ can be expressed by the multiplication $\mu$ and the bullet product $\bullet$. Namely, for
  $f\in C^{m,p}(A,A)$ and $g\in C^{n,q}(A, A)$,
  \begin{equation*}
    f\cup g=(-1)^q(\mu\bullet_{-p} f)\bullet_{m+1} g.
  \end{equation*}
  Indeed, we have
  \begin{equation*}
    \begin{split}
      (\mu\bullet_{-p} f)\bullet_{m+1} g(a_{1, m+n})=& d(\mu\bullet_{-p} f\otimes \id)(a_{1,m}\otimes g(a_{m+1, m+n}))\otimes 1\\
      =&(-1)^pd(\id_p\otimes \mu)(f(a_{1,m})\otimes g(a_{m+1,m+n})\otimes 1\\
      =& (-1)^q f(a_{1,m})g(a_{m+1,m+n}).
    \end{split}
  \end{equation*}
  The differential $\delta$ in $C^{*,*}(A, A)$ can also be expressed using the multiplication and the bullet product $\bullet$. Namely, let $f\in C^{m,p}(A,A)$, then
  \begin{equation*}
    \delta(f)=[f, -\mu]=(-1)^{m-p-1}[\mu, f].
  \end{equation*}
  Indeed, by definition, we have,
  \begin{equation*}
    \begin{split}
    [f,-\mu](a_{1,m+1})=&\sum_{i=1}^{m}(-1)^{i+p}f\bullet_i \mu
    (a_{1,m+1})+
    \sum_{i=1}^2(-1)^{i(m-p-1)+p}\mu\bullet_if(a_{1,m+1})
    +\\
    &\sum_{i=1}^{p}(-1)^{i+m-1}\mu\bullet_{-i}f\\
    =&\sum_{i=1}^m(-1)^{i}f(a_{1,i-1}\otimes a_ia_{i+1}\otimes a_{i+2,m+1})+
     (-1)^{m-1}d(\mu\otimes \id)(f(a_{1,m})\otimes a_{m+1})\otimes 1+\\
     &a_1f(a_{2,m+1})+ \sum_{i=1}^p(-1)^{i+m-1}
     d(\id_i\otimes \mu)(f(a_{1,m})\otimes a_{m+1})\otimes 1  \\
     =& \sum_{i=1}^m(-1)^{i}f(a_{1,i-1}\otimes a_ia_{i+1}\otimes a_{i+2,m+1})+
     a_1f(a_{2,m+1})+ (-1)^mf(a_{1,m})a_{m+1}\\
     =& \delta(f)(a_{1,m})
        \end{split}
  \end{equation*}
  where we used the fact that $d(\overline{\Omega}^p(A))=0$ and $d\circ d=0$ in the third identity. Note that the bullet product $\bullet$ does not define a preLie algebra structure (defined in \cite{Ger1}) on $\Tot(C^*(A, A))$, in general. However, in the following
  proposition we will show that the bullet product $\bullet$ defines a (graded) Lie-admissible algebra structure (defined in Section 2.2 of \cite{MeVa}) on $\Tot(C^*(A, A)).$ That is, the associated Lie bracket $[\cdot,\cdot]$ defines a differential
  graded Lie algebra structure on $\Tot(C^*(A, A))$.
\end{rem}

\begin{prop}\label{prop-lie}
  The bracket $[\cdot,\cdot]$ (defined in  (\ref{equ-defn-bracket})) gives a differential graded Lie algebra (DGLA) structure on the total complex (shifted by $[1]$)
  $$ \Tot(C^{*,*}(A,A))^*[1].$$
   As a consequence,
  \begin{equation*}
  \HH^{>0}(A, \overline{\Omega}^*(A))[1]:=\bigoplus_{\substack{m\in\Z_{>0},\\ p\in Z_{\geq 0}}} \HH^{m}(A, \overline{\Omega}^p(A))[1]
  \end{equation*}
  is a $\Z$-graded Lie algebra, with the grading
  \begin{equation*}
   \HH^{>0}(A, \overline{\Omega}^*(A))_n:=\bigoplus_{\substack{m\in\Z_{>0}, p\in Z_{\geq 0}\\ m-p=n}} \HH^{m}(A, \overline{\Omega}^p(A)).
   \end{equation*}

   \end{prop}
\pf
From the definition of the bracket in  (\ref{equ-defn-bracket}), we observe that
 $[\cdot,\cdot]$ is skew symmetric,
 \begin{equation*}
  [f_1, f_2]=-(-1)^{(m_1-p_1-1)(m_2-p_2-1)}[f_2, f_1].
\end{equation*}
Now let us check the Jacobi identity. That is,
let $f_i\in C^{m_i}(A, \overline{\Omega}^{p_i}(A))$ for $i=1, 2, 3,$
we need to check that
\begin{equation}\label{equ-Jacobi}
  \begin{split}
    (-1)^{n_1n_3}[[f_1, f_2], f_3]+(-1)^{n_2n_1}[[f_2, f_3],
    f_1]+
    (-1)^{n_3n_2}[[f_3, f_1], f_2]=0
  \end{split}
\end{equation}
where $n_i:=m_i-p_i-1,$
for $i=1, 2, 3.$ Now let us apply Lemma \ref{lemma-bullet1} and \ref{lemma-bullet2}, from these two lemmas, we note that
every term on the left hand side of Identity (\ref{equ-Jacobi}) will appear exactly twice.
So the only thing that we should do is to compare the coefficients of the same two terms. However, this
can be done case by case.
  Let us first consider the coefficient of the term $(f_1\bullet_j f_2)\bullet_if_3$ in (\ref{equ-Jacobi}), which is $$(-1)^{n_1n_3+r(m_1,p_1; m_2, p_2;j)+r(m_1+m_2-1, p_1+p_2; m_3,p_3; i)}=(-1)^{n_1n_3+(j-1)n_2+(i-1)n_3+p_3}.$$
By Lemma \ref{lemma-bullet1}, we have the following cases
\begin{enumerate}
\item If $0\leq i-j<m_2$,
we have
\begin{equation}\label{equ-b1}
(f_1\bullet_jf_2)\bullet_i f_3=(-1)^{p_1+p_3}f_1\bullet_j(f_2\bullet_{i-j+1}f_3).
\end{equation}
The coefficient of the term $f_1\bullet_j(f_2\bullet_{i-j+1}f_3)$ is
$$(-1)^{n_1n_3+1+r(m_2,p_2;m_3;p_3;i-j+1)+r(m_1,p_1;m_2+m_3-1;p_2+p_3; j)}=(-1)^{n_1n_3+1+(i-j)n_3+(j-1)(n_2+n_3)+p_1}.$$
Hence the coefficients of these two terms in Identity (\ref{equ-Jacobi}) are up to the scale $-(-1)^{p_1+p_3}$,
then from (\ref{equ-b1}), it follows that these two terms will be cancelled in (\ref{equ-Jacobi}).
\item If $m_2\leq i-j, i<m_1+m_2-p_2,$ then
\begin{equation*}
(f_1\bullet_jf_2)\bullet_if_3=(-1)^{p_3+p_2}(f_1\bullet_{i+p_2-m_2+1} f_3)\bullet_j f_2.
\end{equation*}
The coefficient of the term $(f_1\bullet_{i+p_2-m_2+1} f_3)\bullet_j f_2$ is
$$(-1)^{n_2n_3+1+n_1n_3+r(m_1,p_1;m_3,p_3; i-n_2)+r(m_1+m_3-1,p_1+p_3; m_2,p_2; j)}=(-1)^{n_1n_3+1+
(j-1)n_2+(i-1)n_3+p_2}.$$
Hence the coefficients of these two terms are up to the scale $-(-1)^{p_2+p_3},$ so they will be cancelled in
Identity (\ref{equ-Jacobi}).
\item If $m_2\leq i-j, m_1+m_2-p_2\leq i$,
then  $$(f_1\bullet_jf_2)\bullet_i f_3=f_3\bullet_{-(p_1+p_2+1+i-m_1-m_2)}(f_1\bullet_j f_2).
$$
The coefficient of the term $f_3\bullet_{-(p_1+p_2+1+i-m_1-m_2)}(f_1\bullet_j f_2)$ in (\ref{equ-Jacobi})
is
$$(-1)^{n_2n_3+1+r(m_1,p_1;m_2,p_2;j)+s(m_3,p_3; m_1+m_2-1,p_1+p_2; n_1+n_2-1+i)}=
(-1)^{n_1n_3+1+(j-1)n_2+(i-1)n_3+p_3}.$$
Hence the coefficients are up to the scale $-1$, so in the same reason they will be cancelled.
\item If $1\leq j-i\leq p_3$,
then
$$(f_1\bullet_jf_2)\bullet_i f_3=(-1)^{p_1+p_3}f_1\bullet_i (f_2\bullet_{-(j-i)} f_3).$$
The coefficient of the term $f_1\bullet_i (f_2\bullet_{-(j-i)} f_3)$ in (\ref{equ-Jacobi}) is
$$(-1)^{n_1n_2+1+n_1(n_2+n_3)+r(m_1,p_1;m_2+m_3-1,p_2+p_3; i)+s(m_2,p_2;m_3,p_3; j-i)}=(-1)^{n_1n_3+1+n_2
(i-1)+n_3(j-1)+p_1}.$$
Hence the coefficients are up to the scale $-(-1)^{p_1+p_3}$, so they will be cancelled.
\item If $p_3<j-i$, then
$$(f_1\bullet_jf_2)\bullet_if_3= (-1)^{p_2+p_3}(f_1\bullet_if_3)\bullet_{m_3+j-p_3-1} f_2 $$
The coefficient of the term $(f_1\bullet_if_3)\bullet_{m_3+j-p_3-1} f_2$ is
$$(-1)^{n_2n_3+1+n_1n_3+r(m_1,p_1:m_3,p_3; i)+r(m_1+m_3-1,p_1+p_3; m_2,p_2; n_3+j)}=(-1)^{1+n_1n_3
+n_3(i-1)+(j-1)n_2+p_2}.$$
Hence the coefficients are up to the scale $-(-1)^{p_2+p_3}$, so they will be cancelled.
\end{enumerate}
In a similar way, the other cases can be checked.
So Identity (\ref{equ-Jacobi}) holds.
Hence, it remains to
 verify the following identity, for $f\in C^{m,p}(A,A)$ and $g\in C^{n,q}(A, A)$,
\begin{equation}\label{equ-diff}
\delta([f, g])=(-1)^{n-q-1}[\delta(f), g]+[f, \delta(g)].
\end{equation}
Now by Remark \ref{rem-delta}, it is equivalent to verify that
\begin{equation*}
  [[f,g],\mu]=(-1)^{n-q-1}[[f,\mu], g]+[f, [g,\mu]],
\end{equation*}
which is exactly followed from the Jacobi identity.
Hence Identity (\ref{equ-diff}) holds.
Therefore, $\Tot(C^{*,*}(A, A))^*[1]$ is a differential graded
Lie algebra.
Observe that for any $n\in\Z$, we have,
\begin{equation*}
  H^n(\Tot(C^{*,*}(A, A))^*)\cong \bigoplus_{\substack{m\in\Z_{>0}, p\in Z_{\geq 0}\\ m-p=n}} \HH^{m}(A, \overline{\Omega}^p(A))
\end{equation*}
hence
\begin{equation*}
  \HH^{>0}(A, \overline{\Omega}^*(A))[1]:=\bigoplus_{\substack{m\in\Z_{>0},\\ p\in Z_{\geq 0}}} \HH^{m}(A, \overline{\Omega}^p(A))[1]
\end{equation*} is a $\Z$-graded Lie algebra since the homology of
a differential graded Lie algebra is a graded Lie algebra with the
induced Lie bracket.
\epf

Let $m\in\Z_{>0}$ and $p\in\Z_{\geq 0}$. Then we have
a short exact sequence (from the normalized bar resolution of $A$),
\begin{equation*}
\xymatrix{
0 \ar[r]&  \overline{\Omega}^{p+1}(A)\ar[r]^{\iota_{p+1}}& A^{\otimes (p+2)} \ar[r]^{d_p} & \overline{\Omega}^p(A)\ar[r]& 0,
}
\end{equation*}
which induces a long exact sequence,
\begin{equation}\label{long-new-1}
\xymatrix{
\cdots\ar[r] & \HH^{m}(A, \overline{\Omega}^p(A))\ar[r]^-{\theta_{m,p}} & \HH^{m+1}(A, \overline{\Omega}^{p+1}(A))\ar[r] & \HH^{m+1}(A, A^{\otimes (p+2)}) \ar[r] & \cdots
 }
\end{equation}
where $$\theta_{m,p}:\HH^{m}(A, \overline{\Omega}^p(A))\rightarrow  \HH^{m+1}(A, \overline{\Omega}^{p+1}(A))$$
is the connecting homomorphism. In fact, we can write $\theta_{m,p}$ explicitly.
For any $f\in C^{m}(A, \overline{\Omega}^p(A))$,
\begin{equation}\label{equ-theta-formular}
  \theta_{m,p}(f)(a_{1,m+1})=(-1)^pd(f(a_{1,m})\otimes a_{m+1}\otimes 1).
\end{equation}
Indeed, $\theta_{m, p}$ is induced from the following lifting,
\begin{equation*}
  \xymatrix{
  A\otimes \overline{A}^{\otimes m}\otimes A \ar[d]^{\overline{f}} \ar[dr]^-{f_0}& A\otimes \overline{A}^{\otimes m+1}\otimes A\ar[l]_-{d}\ar[dr]^-{f_1}\\
  \overline{\Omega}^p(A) & A\otimes \overline{A}^{\otimes p}\otimes A\ar@{->>}[l]^-{d} & A\otimes \overline{A}^{\otimes p+1}\otimes A\ar[l]^-{d}
  }
\end{equation*}
where
\begin{equation*}
  \begin{split}
    \overline{f}(a_{1, m+2})&=a_1f(a_{2, m+1})a_{m+2},\\
    f_0(a_{1, m+2})&=(-1)^p a_1f(a_{2, m+1}) \otimes a_{m+2},\\
    f_1(a_{1, m+3})&=(-1)^ma_1f(a_{2, m+1})\otimes a_{m+1}\otimes a_{m+2}.
  \end{split}
\end{equation*}
Hence $\theta_{m, p}(f)(a_{1, m+1})=(-1)^p d(f(a_{1, m})\otimes a_{m+1} \otimes 1).$
As a result, these connecting homomorphisms $\theta_{m,p}$ induce
a homomorphism of degree zero between $\Z$-graded vector spaces,
\begin{equation}\label{equ-connecting}
  \theta: \HH^{*}(A, \overline{\Omega}^*(A))\rightarrow  \HH^{*}(A, \overline{\Omega}^*(A))
\end{equation}
where
\begin{equation*}
  \theta|_{\HH^{m}(A, \overline{\Omega}^p(A))}=\theta_{m, p}:\HH^{m}(A, \overline{\Omega}^p(A))\rightarrow  \HH^{m+1}(A, \overline{\Omega}^{p+1}(A)).
\end{equation*}
The following proposition shows that $\theta$ is a  module homomorphism of the $\Z$-graded
Lie algebra $\HH^*(A,\overline{\Omega}^*(A))$.
\begin{prop}\label{prop-hom-ger}
Let $A$ be an associative algebra over a commutative algebra $k$, then the homomorphism of $\Z$-graded $k$-modules
  \begin{equation*}
  \theta: \HH^{*}(A, \overline{\Omega}^*(A))\rightarrow  \HH^{*}(A, \overline{\Omega}^*(A))
\end{equation*}
(defined in (\ref{equ-connecting}) above) is a module homomorphism of degree zero over the
$\Z$-graded Lie algebra $\HH^{*}(A, \overline{\Omega}^*(A))$.
\end{prop}
\pf Let $f\in\HH^m(A, \overline{\Omega}^p(A))$ and $g\in\HH^n(A,\overline{\Omega}^q
(A))$, it is sufficient to verify that
\begin{equation}\label{equ-77}
  \theta_{m+n-1,p+q}([f, g])=[\theta_{m,p}(f),g].
\end{equation}
First we claim that the following two identities hold,
\begin{equation}\label{equ-78}
\begin{split}
\theta_{m,p}(f)\bullet g-\theta_{m+n-1,p+q}(f\bullet g)&=(-1)^{m(q-n-1)+p+q}
d(f\otimes  g\otimes 1)\\
g\bullet \theta_{m,p}(f)-\theta_{m+n-1, p+q}(g\bullet f)&=(-1)^{(p+1)(q-n-1)+p+q} d(f\otimes g\otimes 1).
 \end{split}
\end{equation}
It is easy to check that (\ref{equ-78}) implies Identity (\ref{equ-77}). Now let us
verify the two identities in (\ref{equ-78}). Indeed, we have
\begin{equation*}
  \begin{split}
    \theta_{m,p}(f)\bullet g=&\sum_{i=1}^{m+1}(-1)^{r(m+1, p+1;n, q;i)}
    \theta_{m,p}(f)\bullet_i g+\sum_{i=1}^q(-1)^{s(m+1,p+1;n, q;i)} \theta_{m,p}(f)\bullet_{-i} g\\
    =&\sum_{i=1}^{m}(-1)^{r(m, p;n, q;i)}
    \theta_{m+n-1,p+q}(f\bullet_i g)+\sum_{i=1}^q(-1)^{s(m,p;n, q;i)} \theta_{m,p}(f)\bullet_{-i} g+\\
    &(-1)^{r(m+1,p+1;n,q;m+1)}\theta_{m,p}(f)\bullet_{m+1} g\\
    =& \theta_{m+n-1,p+q}(f\bullet g)+(-1)^{m(q-n-1)+p+q}d(f\otimes g\otimes 1)
  \end{split}
\end{equation*}
Similarly, we have
\begin{equation*}
  \begin{split}
    g\bullet \theta_{m,p}(f)=&
   \sum_{i=1}^n (-1)^{r(n,q;m+1,p+1;i)}g\bullet_i\theta_{m,p}(f)+\sum_{i=1}
   ^{p+1}(-1)^{s(n,q;m+1,p+1;i)}g\bullet_{-i}\theta_{m,p}(f)\\
   =&\sum_{i=1}^n(-1)^{r(n,q;m,p;i)}\theta_{m+n-1,p+1}(g\bullet_i f)+\sum_{i=1}^p
   (-1)^{s(n,q;m,p;i)}\theta_{m+n-1,p+1}(g\bullet_{-i} f)+\\
   &(-1)^{s(n,q;m+1,p+1;p+1)}g\bullet_{-p-1}\theta_{m,p}(f)\\
   =&\theta_{m+n-1,p+q}(g\bullet f)+(-1)^{(p+1)(q-n-1)+p+q}d(f\otimes
   g\otimes 1).
     \end{split}
\end{equation*}
Hence we have proved that the two identities in (\ref{equ-78}) hold, so the proof is completed.
\epf
\begin{rem}\label{rem-H0}
  Note that we did not consider $\HH^0(A, \Omega^p(A))$ for $p\in\Z_{>0}$ when we defined the Lie bracket $[\cdot,\cdot]$ on $\HH^{>0}(A, \overline{\Omega}^*(A))$.
  In fact, for $p\in\Z_{\geq 0}$, let us consider the homomorphism
  $$\theta_{0,p}: \HH^0(A, \overline{\Omega}^p(A))\rightarrow \HH^1(A, \overline{\Omega}^{p+1}(A)).$$
  Via this homomorphism, we can define a Lie bracket action of $\HH^0(A, \Omega^p(A))$ on $\HH^n(A, \Omega^q(A))$ for $n\in\Z_{>0}$ and $p\in\Z_{\geq 0}$ as follows:
  for $\alpha\in  \HH^0(A, \Omega^p(A))$ and $g\in \HH^n(A, \Omega^q(A)),$ define
  \begin{equation*}
    [\alpha, g]:=[\theta_{0,p}(\alpha), g].
  \end{equation*}
Denote
\begin{equation*}
   \HH^{\geq 0}(A, \overline{\Omega}^*(A)):=\bigoplus_{\substack{m\in\Z_{\geq 0}, p\in Z_{\geq 0}}} \HH^{m}(A, \overline{\Omega}^p(A)).
   \end{equation*}
In the following proposition, we will prove that there is a Gerstenhaber algebra structure
on $\HH^{\geq 0}(A, \overline{\Omega}^*(A))$.
\end{rem}

\begin{prop}\label{prop-ger3}
  Let $A$ be an associative algebra over a commutative ring $k$. Suppose that $A$ is projective as a $k$-module.
  Then
  \begin{equation*}
   \HH^{\geq 0}(A, \overline{\Omega}^*(A)):=\bigoplus_{\substack{m\in\Z_{\geq 0}, p\in Z_{\geq 0}}} \HH^{m}(A, \overline{\Omega}^p(A))
   \end{equation*}
   with the Lie bracket $[\cdot,\cdot]$ and the cup product $\cup$, is a Gerstenhaber algebra.
\end{prop}
\pf
From Proposition \ref{prop-lie}, Proposition \ref{prop-hom-ger} and Remark \ref{rem-H0}, it follows that $\HH^{\geq 0}(A, \overline{\Omega}^*(A))$ is a
$\Z$-graded Lie algebra. Next let us prove that $\cup$ defines a graded commutative algebra structure.
Let $f_i\in\HH^{m_i}(A, \overline{\Omega}^{p_i}(A))$ for $i=1, 2$. Then we claim that
\begin{equation*}
  \begin{split}
    f_1\cup f_2-(-1)^{(m_1-p_1)(m_2-p_2)} f_2
    \cup f_1
    =\delta(f_2\bullet f_1),
  \end{split}
\end{equation*}
indeed, by immediate calculation,  we have
\begin{equation*}
  \begin{split}
    \delta(f_2\bullet f_1)(a_{1, m+n})=&\sum_{i=1}^{m_2}(-1)^{p_1+p_2+(i-1)(p_1-m_1-1)}\delta(f_2\circ_i f_1)(a_{1, m+n}) +\\
    &\sum_{i=1}^{p_1} (-1)^{p_1+p_2+i(p_2-m_2-1)} \delta(f_2\circ_{-i} f_1)(a_{1, m+n})\\
   =& f_1(a_{1, m_1})f_2(a_{m_1+1, m_1+m_2})-(-1)^{(m_1-p_1)(m_2-p_2)} f_2(a_{1, m_2})f_1(a_{m_2+1, m_1+m_2}).
  \end{split}
\end{equation*}
Hence we obtain, in $\HH^{m_1+m_2}(A, \overline{\Omega}^{p_1+p_2}(A))$
$$f_1\cup f_2=(-1)^{(m_1-p_1)(m_2-p_2)} f_2\cup f_1.$$
So the graded commutativity of cup products holds.
It remains to verify the compatibility between $\cup$ and $[\cdot,\cdot]$, namely, for $f_i\in\HH^{m_i}(A,\overline{\Omega}^{p_i}(A))$,
\begin{equation}\label{equ-compa}
  [f_1\cup f_2, f_3]=[f_1, f_3]\cup f_2+(-1)^{(m_3-p_3-1)(m_1-p_1)}f_1\cup [f_2, f_3].
\end{equation}
\begin{claim}\label{claim-equ2}
In the cohomology group $\HH^{m_1+m_2+m_3-1}(A,\overline{\Omega}^{p_1+p_2+p_3}(A))$, we also have the following
identity,
\begin{equation}\label{equ-claim2}
f_3\bullet (f_1\cup f_2)-(-1)^{(m_2-p_2)(m_3-p_3-1)}(f_3\bullet f_1)\cup f_2-f_1\cup(f_3\bullet f_2)=0.
\end{equation}
\end{claim}
\begin{claim}\label{claim-equ1}
In the cohomology group $\HH^{m_1+m_2+m_3-1}(A,\overline{\Omega}^{p_1+p_2+p_3}(A))$, we have the following
identity,
\begin{equation}\label{equ-claim1}
(f_1\cup f_2)\bullet f_3-(f_1\bullet f_3)\cup f_2-(-1)^{(m_1-p_1)(m_3-p_3-1)}f_1\cup (f_2\bullet f_3)=0.
\end{equation}
\end{claim}
It is easy to check that these two claims imply Identity (\ref{equ-compa}).
Now let us prove these two claims. The proofs of these two claims are very similar to
the proof of Theorem 5 in \cite{Ger1}.

{\it Proof of Claim \ref{claim-equ2}, \hspace{2ex}}
First, it is easy to check that we have the following identity for the left hand side
in (\ref{equ-claim2})
\begin{equation*}
  \begin{split}
    \mbox{LHS}=&\sum_{i=1}^{m_3} (-1)^{r(m_3,p_3;m_1+m_2,p_1+p_2;i)}f_3\bullet_i(f_1\cup f_2)+
\sum_{i=1}^{p_1} (-1)^{s(m_3,p_3;m_1+m_2,p_1+p_2;i)}f_3\bullet_{-i}(f_1\cup f_2)-\\
    & \sum_{i=1}^{m_3}(-1)^{(m_2-p_2)(m_3-p_3-1)+r(m_3,p_3;m_1,p_1;
    i)}(f_3\bullet_if_1)\cup f_2- \\  &\sum_{i=1}^{p_1}(-1)^{(m_2-p_2)(m_3-p_3-1)+s(m_3,p_3;m_1,p_1;
    i)}(f_3\bullet_{-i} f_1)\cup f_2-
     \sum_{i=1}^{m_3} (-1)^{r(m_3,p_3;m_2,p_2;i)}f_1\cup (f_3\bullet_i
    f_2).
  \end{split}
\end{equation*}
Set
\begin{equation*}
  \begin{split}
    H_1:=& \sum_{i=1}^{p_1}\sum_{j=1}^{m_3-1}(-1)^{\epsilon_{i,j}} d(\id^{\otimes i}\otimes f_3\otimes \id^{\otimes p_1+p_2-i})(f_1\otimes \id^{\otimes j-1} \otimes f_2\otimes \id^{\otimes m_3-j-1})\otimes 1+\\
   & \sum_{i=1}^{m_3}\sum_{j=1}^{m_2-1} (-1)^{\epsilon'_{i,j}}d((f_3\otimes \id^{\otimes p_1+p_2})(\id^{\otimes i-1}\otimes f_1\otimes \id^{\otimes j-1} \otimes f_2\otimes \id^{\otimes m_3-i-j})\otimes 1).
  \end{split}
\end{equation*}
where to simplify the formula, we do not consider the sign of each term in
$H_1$.
Then we will show that
$$\mbox{LHS}=\delta(H_1).$$
Namely, take any element $a_{1, m_1+m_2+m_3-1}\in A^{\otimes m_1+m_2+m_3-1}$,
we need to prove that 
\begin{equation}\label{equ-com}
  \mbox{LHS}(a_{1, m_1+m_2+m_3-1})=\delta(H_1)(a_{1, m_1+m_2+m_3-1}).
\end{equation}
Indeed, this identity can be proved by a recursive procedure. 
That is, as a first step, we will verify that those terms containing 
$f_1(a_{1, m_1})$ in (\ref{equ-com}) can be cancelled by frequently using
the fact that $d^2=0$ and $d(\overline{\Omega}^p(A))=0$ for $p\in \Z_{>0}$.
This can be done by comparing the terms containing $f_1(a_{1, m_1})$ of 
both sides in (\ref{equ-com}).
Then after cancelling those terms containing $f_1(a_{1, m_1})$,
we obtain a new identity $$\mbox{LHS}'(a_{1, m_1+m_2+m_3-1})=\delta(H_1)'(a_{1, m_1+m_2+m_3-1}),$$
which does not have the terms containing $f_1(a_{1, m_1})$. By similar processing, 
we will cancel the terms containing $f_1(a_{2, m_1+1})$. Hence after several times,
all the terms will be cancelled, so Identity (\ref{equ-com}) holds.
Therefore $\mbox{LHS}=0$ in the cohomology group $\HH^{m_1+m_2+m_3-1}(A, \overline{\Omega}^{p_1+p_2+p_3}(A))$.
\epf

{\it Proof of Claim \ref{claim-equ1}, \hspace{2ex}} Similarly,  the left hand side in (\ref{equ-claim1}) can be written as
\begin{equation*}
\begin{split}
\mbox{LHS}=&\sum_{i=1}^{m_1}(-1)^{r(m_1+m_2,p_1+p_2;m_3,p_3;i)}
(f_1\cup f_2)\bullet_i f_3+\sum_{i=1}^{p_3}(-1)^{s(m_1+m_2,p_1+p_2;m_3,p_3;i)}
(f_1\cup f_2)\bullet_{-i} f_3-\\
&\sum_{i=1}^{m_1} (-1)^{r(m_1,p_1;m_3,p_3;i)}(f_1\bullet_if_3)\cup f_2-\sum_{i=1}^{p_3} (-1)^{s(m_1,p_1;m_3,p_3;i)}(f_1\bullet_{-i}f_3)\cup f_2-\\
&\sum_{i=1}^{p_3}(-1)^{(m_1-p_1)(m_3-p_3-1)+s(m_2,p_2;m_3,p_3; i)}f_1\cup(f_2\bullet_{-i} f_3).
\end{split}
\end{equation*}
Set
\begin{equation*}
  \begin{split}
    H_2:=&\sum_{i=1}^{m_1}\sum_{j=1}^{p_3} (-1)^{\epsilon_{i, j}} d((f_1\otimes\id^{\otimes j-1}\otimes f_2\otimes \id^{\otimes p_3-j})(\id^{\otimes i-1}\otimes f_3\otimes \id^{\otimes m_1+m_2-i-1})\otimes 1)+\\
    &\sum_{i=1}^{p_3}(-1)^{\epsilon_i}d((\id\otimes f_1\otimes \id^{\otimes i-1} \otimes f_2\otimes \id^{\otimes p_3-i-1})(f_3(a_{1, m_3})\otimes a_{m_3, m_1+m_2+m_3-2})\otimes 1)
  \end{split}
\end{equation*}
Similarly as above, by calculation, we obtain that
$$\mbox{LHS}=\delta(H_2),$$
hence $\mbox{LHS}=0$ in the cohomology group $\HH^{m_1+m_2+m_3-1}(A, \overline{\Omega}^{p_1+p_2+p_3}(A))$.
\epf\\
Therefore, we have completed the proof.
\epf
\begin{rem}
  The proof of Proposition \ref{prop-ger3} above relies on combinatorial calculations. To understand the Gerstenhaber algebra structure much better, it is interesting to investigate whether
  there is a $B_{\infty}$-algebra structure on the total complex of  $C^{*, *}(A, A)$ since from
  Section 5.2 in \cite{GeJo} it follows that a $B_{\infty}$-algebra structure on a chain complex $C$ induces
  a canonical Gerstenhaber algebra structure on the homology $H_*(C)$.

\end{rem}
Now let us prove the main theorem (cf. Theorem \ref{thm-gerst}) in this section.

{\it Proof of Theorem \ref{thm-gerst}, \hspace{2ex}}
From Proposition \ref{prop}, we have the following isomorphism for $m\in \Z$,
$$\Phi_m: \lim_{\substack{\longrightarrow\\p\in\Z_{\geq 0}\\ m+p>0}} \HH^{m+p}(A, \overline{\Omega}^p(A))\rightarrow \Hom_{\DD_{\sg}(A^e)}(A, A[m]).$$
From Proposition \ref{prop-hom-ger} and Proposition \ref{prop-ger3}, it follows that the structural morphism in the direct system $\HH^{m+*}(A, \overline{\Omega}^{p+*}(A)),$
$$\theta_{m, p}:\HH^{m}(A, \overline{\Omega}^p(A))\rightarrow  \HH^{m+1}(A, \overline{\Omega}^{p+1}(A))$$
preserves the Gerstenhaber algebra structure. Therefore, there is an induced Gerstenhaber algebra structure on its direct limit $$\HH^*_{\sg}(A, A)\cong \lim_{\substack{\longrightarrow\\p\in\Z_{\geq 0}\\ *+p>0}} \HH^{*+p}(A, \overline{\Omega}^p(A)).$$
So $\HH_{\sg}^*(A, A)$ is a Gerstenhaber algebra, equipped with the cup product
$\cup$ and the induced Lie bracket $[\cdot, \cdot]$.
\epf

Let us denote, for $m\in \Z_{\geq 0},$
$$\Ker^{m, \infty}(A, A):=\ker(\HH^m(A, A)\rightarrow \HH^m_{\sg}(A, A)),$$
and
$$\HH^{m, \infty}(A, A):=\Imm(\HH^m(A, A)\rightarrow \HH^m_{\sg}(A, A)).$$
Then we have
$$\HH^{m, \infty}(A, A)=\HH^m(A, A)/\Ker^{m, \infty}(A, A).$$
\begin{cor}
  $\Ker^{*, \infty}(A, A)$ is a Gerstenhaber ideal of $\HH^*(A, A)$. In particular,
  $\HH^{*, \infty}(A, A)$ is a Gerstenhaber subalgebra of $\HH^*_{\sg}(A, A)$.
\end{cor}
\pf First let us claim that the natural morphism
$$\Phi: \HH^*(A, A)\rightarrow \HH_{\sg}^*(A, A)$$
is a Gerstenhaber algebra homomorphism.  In order to prove this claim, from Proposition \ref{prop-hom-ger} and Remark \ref{rem-H0},  it is sufficient to
verify that
\begin{equation*}
  [\theta_{0, 0}(\HH^0(A, A)), \theta_{0, 0}(\HH^0(A, A))]=0,
\end{equation*}
where we recall that $$\theta_{0, 0}: \HH^0(A, A)\rightarrow \HH^1(A, \overline{\Omega}^1(A))$$
is the connecting homomorphism in (\ref{long-new-1}).
Let $\lambda,\mu\in \HH^0(A, A)$, then from (\ref{equ-theta-formular}) it follows that for any $a\in A$,
\begin{equation*}
\begin{split}
  \theta_{0, 0}(\lambda)(a)=\lambda\otimes a -\lambda a\otimes 1;\\
  \theta_{0, 0}(\mu)(a)=\mu\otimes a -\mu a\otimes 1.\\
\end{split}
\end{equation*}
Hence, by direct calculation,  we obtain that for any $a\in A$,
\begin{equation*}
  \begin{split}
    [\theta_{0, 0}(\lambda),\theta_{0, 0}(\mu)](a)=0.
  \end{split}
\end{equation*}
So we have shown that
$$\Phi: \HH^*(A, A)\rightarrow \HH_{\sg}^*(A, A)$$
is a Gerstenhaber algebra homomorphism. Thus its kernel $\Ker^{*, \infty}(A, A)$ is a Gerstenhaber ideal and its image $\HH^{* \infty}(A, A)$
is a Gerstenhaber algebra.
\epf

\section{Gerstenhaber algebra structure on $\HH^*(A, A^{\otimes >0})$}

Let $m,n\in\Z_{>0}$ and $p, q\in \Z_{>1}$, we will define a star product as follows,
\begin{equation*}
  \star:C^m(A, A^{\otimes p})\times C^n(A, A^{\otimes q})\rightarrow C^{m
  +n-1}(A, A^{\otimes p+q-1})
\end{equation*}
for $f\in C^m(A, A^{\otimes p})$ and $g\in C^n(A, A^{\otimes q})$,
denote
\begin{equation*}
\begin{split}
  f\star_0 g:&=(f\otimes \id)(\id_{m-1}\otimes g),\\
  f\star_1 g:&=(\id_{p-1}\otimes f)(g\otimes \id).
\end{split}
\end{equation*}
Then we define
\begin{equation*}
\begin{split}
  f\star g:=(-1)^{(m-1)(n-1)}f\star_0 g+f\star_1 g
\end{split}
\end{equation*}
%
and denote
\begin{equation}\label{equ-bra-a}
  \{f, g\}:=f\star g-(-1)^{(m-1)(n-1)}g\star f.
\end{equation}

\begin{rem}
For the case $p=1$, we can also define a Lie bracket
 \begin{equation}\label{equ-bra-b}
   \{\cdot,\cdot\}: C^m(A, A)\times C^n(A, A^{\otimes q})\rightarrow C^{m+n-1}(A, A^{\otimes q}),
  \end{equation}
  let $f\in C^m(A, A)$ and $g\in C^n(A, A^{\otimes q})$,
  \begin{equation*}
    \{f, g\}=f\star g-(-1)^{(m-1)(n-1)}g\circ f,
  \end{equation*}
  where $g\circ f$ is the circ product (cf. (\ref{equ-circ-product})) defined in \cite{Ger1}.
\end{rem}
We remark that in general, the star product $\star$ is not associative.
However, some associativity properties hold. In the following lemma, we list
some of them.

\begin{lemma}\label{lemma-star}
 Let $f_i\in C^{m_i}(A,
A^{\otimes p_i}), i=1, 2, 3.$  Then we have the following,
\begin{enumerate}
  \item for $p_1\in\Z_{>0}, p_2,p_3\in\Z_{>1}$,
   \begin{equation*}
   \begin{split}
   (f_1\star_0 f_2)\star_0 f_3=f_1\star_0(f_2\star_0 f_3),\\
   (f_1\star_0 f_2)\star_1 f_3=(f_1\star_1 f_3)\star_0 f_2,\\
   (f_1\star_1 f_2)\star_1 f_3=f_1\star_1(f_2\star_1 f_3),
   \end{split}
  \end{equation*}
  \item for $p_1,p_2\in\Z_{>0}, p_3\in\Z_{>1}$,
   $$f_1\star_1(f_2\star_0 f_3)=f_2\star_0(f_1\star_1 f_3).$$
\end{enumerate}
\end{lemma}

Note that $C^{*>0}(A, A^{\otimes *> 0})$  is a double complex with horizontal and vertical differentials induced from the bar resolution. We consider the total complex
\begin{equation*}
 \Tot(C^{*>0}(A, A^{\otimes *> 0}))^*,
\end{equation*}
which is a differential $\Z$-graded $k$-module, with the grading
\begin{equation*}
  \Tot(C^{*>0}(A, A^{\otimes *>1}))^m=\bigoplus_{\substack{p\in\Z_{>0}}}C^m(A, A^{\otimes p})
\end{equation*}
and the differential induced by the horizontal differential.
\begin{rem}\label{rem-prod-delta}
We remark that the horizontal differential $\delta$ in $C^{*>0}(A, A^{*>0})$ can be expressed by the
star product, circ product and the multiplication $\mu$ of $A$. That is, for $f\in C^m(A, A^{\otimes p})$,
\begin{equation*}
\delta(f)=(-1)^{m-1}\{\mu, f\}.
\end{equation*}
Indeed, we have,
\begin{equation*}
\begin{split}
\{\mu, f\}(a_{1, m+1})=&\mu\star f(a_{1, m+1})-(-1)^{m-1}f\circ \mu(a_{1, m+1})\\
=&(-1)^{m-1}(\mu\otimes \id)(a_1\otimes f(a_{2, m+1}))+(\id_{p-1}\otimes \mu)(f(a_{1,m})\otimes a_{m+1})-\\
&(-1)^{m-1}\sum_{i=1}^m(-1)^{i-1}f(a_{1,i-1}\otimes a_ia_{i+1}\otimes a_{i+2, m+1})\\
=&(-1)^{m-1}\delta(f)(a_{1,m+1}).
\end{split}
\end{equation*}
\end{rem}
\begin{prop}
  $\Tot(C^{*>0}(A, A^{\otimes *>0}))^*$ is a differential graded Lie algebra (DGLA) with the bracket $\{\cdot, \cdot\}$ defined in (\ref{equ-bra-a}). As a consequence,
  \begin{equation*}
  \HH^{*>0}(A, A^{\otimes *>0})=\bigoplus_{m,p\in\Z_{>0}}\HH^m(A, A^{\otimes p})
  \end{equation*}
  is a $\Z$-graded Lie algebra with the grading
  \begin{equation*}
  \HH^{*>0}(A, A^{\otimes *>0})^m:=\bigoplus_{p\in\Z_{>0}} \HH^m(A, A^{\otimes p}).
    \end{equation*}
\end{prop}
\pf
Let us prove first that it is a $\Z$-graded Lie algebra. Note that the bracket is skew-symmetric. We need to verify the Jacobi identity, namely, for $f_i\in C^{m_i}(A,
A^{\otimes p_i}), i=1, 2, 3,$
\begin{equation}\label{equ-Jacobi2}
  \begin{split}
    (-1)^{n_1n_3}\{\{f_1, f_2\}, f_3\}+(-1)^{n_2n_1}\{\{f_2, f_3\},
    f_1\}+
    (-1)^{n_3n_2}\{\{f_3, f_1\}, f_2\}=0
  \end{split}
\end{equation}
where $n_i:=m_i-1.$
 Let us discuss the following three cases.

{\bf Case 1:} $p_i\in\Z_{>1}$.
From Lemma \ref{lemma-star}, it follows that each term of the right
hand side in (\ref{equ-Jacobi2}) will appear exactly twice, hence it is sufficient to compare the coefficients of those two same terms.
For example,
the coefficient of $(f_1\star_0 f_2)\star_0 f_3$ in (\ref{equ-Jacobi2}) is
$(-1)^{n_2(n_1+n_3)}.$
From Lemma \ref{lemma-star}, we obtain that
\begin{equation*}
  (f_1\star_0 f_2)\star_0 f_3=f_1\star_0 (f_2\star_0 f_3).
\end{equation*}
And the coefficient of $f_1\star_0 (f_2\star_0 f_3)$ in (\ref{equ-Jacobi2}) is
 $ -(-1)^{n_2(n_1+n_3)},$
hence these two same terms will be cancelled in (\ref{equ-Jacobi2}).
Similarly, we can verify the other terms.

{\bf Case 2:} $p_1,p_2=1$ and $p_3\in\Z_{>1}.$
From Theorem 2 in \cite{Ger1}, we have the following identity,
$$f_3\circ\{f_1, f_2\}=(f_3\circ f_1)\circ f_2-(-1)^{(m_1-1)(m_2-1)}(f_3\circ f_2)\circ f_1.$$
Thus, to verify the Jacobi identity (\ref{equ-Jacobi2}) is equivalent to verify the following one,
\begin{equation}\label{99}
\begin{split}
\{f_1, f_2\}\star f_3 =&f_1\star (f_2\star f_3 )-(-1)^{(m_2-1)m_3}f_1\star (f_3 \circ f_2)-(-1)^{(m_1-1)(m_2+m_3-1)}(f_2\star f_3 )\circ f_1-\\
&(-1)^{(m_1-1)(m_2-1)}f_2\star (f_1\star f_3 )+(-1)^{(m_1-1)m_2}f_2\star (f_3 \circ f_1)+\\
&(-1)^{(m_3-1)(m_2-1)}(f_1\star f_3 )\circ f_2.
\end{split}
\end{equation}
It is easy, by definition, to check that the following identity holds,
\begin{equation}\label{100}
  \begin{split}
   & (f_1\circ f_2)\star f_3+(-1)^{(m_2-1)m_3}f_1\star (f_3 \circ f_2)-(-1)^{(m_3-1)(m_2-1)}(f_1\star f_3 )\circ f_2\\
   =&(-1)^{(m_1+m_2)(m_3-1)+(m_1-1)(m_2-1)} f_1\star_0(f_2\star_0f_3)+
   f_1\star_1(f_2\star_1 f_3).
  \end{split}
\end{equation}
Similarly,  we also have
\begin{equation}\label{101}
  \begin{split}
   & (f_2\circ f_1)\star f_3+(-1)^{(m_1-1)m_3}f_2\star (f_3 \circ f_1)-(-1)^{(m_3-1)(m_1-1)}(f_2\star f_3 )\circ f_1\\
  =&(-1)^{(m_1+m_2)(m_3-1)+(m_1-1)(m_2-1)} f_2\star_0(f_1\star_0f_3)+
    f_2\star_1(f_1\star_1 f_3).
  \end{split}
\end{equation}
We also note that
\begin{equation}\label{102}
  \begin{split}
     &(f_1\star (f_2\star f_3)-(-1)^{(m_1-1)(m_2-1)}f_2\star (f_1\star f_3))\\
     =&(-1)^{(m_1-1)(m_2+m_3)}f_1\star_0(f_2\star_0 f_3)+
     f_1\star_1(f_2\star_1 f_3)-\\
     &(-1)^{(m_1-1)(m_2+m_3)+(m_2-1)(m_1+m_3
     )}f_2\star_0(f_1\star_0 f_3)-(-1)^{(m_1-1)(m_2-1)}
     f_2\star_1(f_1\star_1 f_3).
  \end{split}
\end{equation}
So combining (\ref{100}), (\ref{101}) and (\ref{102}), we get
Identity (\ref{99}). Hence the Jacobi identity holds.

{\bf Case 3:} $p_1=1$, $p_2, p_3\in\Z_{>1}$.
First, similar to Case 2, we have the following identities for $\{i, j\}=\{ 2, 3\}$.
\begin{equation}
  \begin{split}
    &(f_i\circ f_1)\star f_j+(-1)^{(m_1-1)m_j}f_2\star(f_3\circ f_1)-(-1)^{(m_1-1)(m_j-1)}
    (f_i\star f_j)\circ f_1\\
    =&(-1)^{(m_1+m_i)(m_j-1)+(m_1-1)(m_i-1)}f_2\star_0(f_1\star_0 f_j)-f_i\star_1(f_1\star_1 f_j).
       \end{split}
\end{equation}
Hence, each term in Jacobi identity in (\ref{equ-Jacobi2}) can be expressed by star product (no circ product).
So, similar to Case 1, by using Lemma \ref{lemma-star}, we can compare the same two terms.
Thus, we have verified that
$\Tot(C^{*>0}(A, A^{\otimes *> 0}))^*$ is a $\Z$-graded Lie algebra.

Next we will prove that $\{\cdot,\cdot\}$ is compatible with differential, that is,
let $f_i\in C^{m_i}(A, A^{\otimes p_i})$ for $i=1, 2,$ we need to verify the following
identity,
\begin{equation}
  \delta(\{f_1, f_2\})=(-1)^{m_1-1}\{f_1, \delta(f_2 )\}+\{\delta(f_1), f_2\}.
\end{equation}
From Remark \ref{rem-prod-delta}, we note that it follows from the Jacobi identity (\ref{equ-Jacobi2}).
Since $$H^m(\Tot(C^{*>0}(A, A^{\otimes *> 0}))\cong \HH^{*>0}(A, A^{\otimes *>0})^m,$$
it follows that
$\HH^{*>0}(A, A^{\otimes *>0})$ is a $\Z$-graded Lie algebra.
Hence
we have completed the proof.
\epf

\begin{lemma}\label{lemma-zero-product}
Let $f_i\in\HH^{m_i}(A, A^{\otimes p_i})$ for $i=1, 2$ and suppose that $p_2>1$. Then we have
$f_1\cup f_2=f_2\cup f_1=0$.
\end{lemma}
\pf In fact, we have the following identity for $f_i\in C^{m_i}(A, A^{\otimes p_i})$ and $p_2>1$,
\begin{equation*}
\begin{split}
\delta(f_1\star_0 f_2)&=\delta(f_1)\star_0 f_2+(-1)^{m_1-1} f_1\star_0 \delta(f_2)+(-1)^{m_1} f_1\cup f_2,\\
\delta(f_1\star_1 f_2)&=f_1\star_1 \delta(f_2)+(-1)^{m_2-1} \delta(f_1)\star_1 f_2+(-1)^{m_2} f_2\cup f_1.
\end{split}
\end{equation*}
Hence it follows that $f_1\cup f_2=f_2\cup f_1=0.$
\epf

Next we will prove that $\HH^{*>0}(A, A^{\otimes *>0})$ is a Gerstenhaber algebra.
\begin{prop}\label{prop-zero-product}
$\HH^{*>0}(A, A^{\otimes *>0})$ is a Gerstenhaber algebra (without the unity) with the cup product and Lie bracket.
\end{prop}
\pf It remains to verify the compatibility between cup product and Lie bracket. On the other hand,
from Lemma \ref{lemma-zero-product} above, it is sufficient to verify that for $f_i\in \HH^{m_i}(A, A)$,
$i=1, 2$ and $f_3\in \HH^{m_3}(A, A^{\otimes p_3})$, where $p_3>1$,
\begin{equation}\label{equ-com-bar}
\begin{split}
\{f_1\cup f_2, f_3\}&=(-1)^{m_1(m_3-1)}f_1\cup\{f_2, f_3\}+(-1)^{(m_1+m_3-1)m_2} f_2\cup\{f_1, f_3\}\\
\end{split}
\end{equation}
Recall that we have the following identity in $\HH^{m_1+m_2+n-1}(A, A^{\otimes p})$(cf. Theorem 5. \cite{Ger1}),
\begin{equation}\label{id2}
\begin{split}
f_3\circ (f_1\cup f_2)-f_1\cup (f_3\circ f_2)-(-1)^{m_2(m_3-1)}(f_3\circ f_1)\cup f_2=0.\\
\end{split}
\end{equation}
Hence it follows that Identity (\ref{equ-com-bar}) is equivalent to the following identity,
\begin{equation}\label{identity14}
\begin{split}
(f_1\cup f_2)\star f_3-(f_1\star f_3)\cup f_2-(-1)^{m_1(m_3-1)}f_1\cup (f_2\star f_3)=0.\\
\end{split}
\end{equation}
Let us compute the left hand side in (\ref{identity14}),
\begin{equation}
  \begin{split}
    &((f_1\cup f_2)\star f_3-(f_1\star f_3)\cup f_2-(-1)^{m_1(n-1)}f_1\cup (f_2\star f_3))\\
   =&(-1)^{(m_1+m_2-1)(m_3-1)}f_1\cup (f_2\star_0 f_3)+(f_1\star_1 f_3)\cup f_2-(-1)^{(m_1-1)(m_3-1)}(f_1\star_0 f_3)\cup f_2-\\
  & (f_1\star_1 f_3)\cup f_2-(-1)^{m_1(m_3-1)+(m_2-1)(m_3-1)} f_1\cup (f_2\star_0 f_3)-(-1)^{m_1(m_3-1)}f_1\cup(f_2\star_1 f_3)\\
  =&-(-1)^{(m_1-1)(m_3-1)}(f_1\star_0 f_3)\cup f_2 -(-1)^{m_1(m_3-1)}f_1\cup(f_2\star_1 f_3)
    \end{split}
\end{equation}
Set
\begin{equation*}
  G:=(-1)^{m_1m_3}(f_1\star_0 f_3)\star_1 f_2.
\end{equation*}
By calculation, we obtain that
\begin{equation*}
  \delta(G)=(-1)^{(m_1-1)(m_3-1)}(f_1\star_0 f_3)\cup f_2 +(-1)^{m_1(m_3-1)}f_1\cup(f_2\star_1 f_3),
\end{equation*}
hence it follows that in $\HH^{m_1+m_2+n-1}(A, A^{\otimes p})$,
\begin{equation*}
  (f_1\cup f_2)\star g-(f_1\star g)\cup f_2-(-1)^{m_1(n-1)}f_1\cup (f_2\star g)=0.
\end{equation*}
Therefore, we have verified Identity (\ref{equ-com-bar}).
\epf

\begin{cor}
Let $A$ be an associative algebra over a commutative ring $k$. Let $f_i\in\HH^{m_i}(A, A)$,
$m_i\in\Z_{>0}$ for $i=1, 2$. Then for any $g\in \HH^m(A, A^{\otimes p}), p\in \Z_{>1}, m>0,$ we have
\begin{equation*}
\{f_1\cup f_2, g\}=0.
\end{equation*}
\end{cor}
\pf This result follows from Lemma \ref{lemma-zero-product} and Proposition \ref{prop-zero-product}.
\epf

%

\section{Special case: self-injecitve algebra}

In this section, let $A$ be a self-injective algebra over a field $k$.
Then $A^e:=A\otimes_kA^{\op}$ is also a self-injective algebra.
Recall that the singular category $\DD_{\sg}(A)$ of a self-injective algebra $A$ is equivalent to the stable module category $A$-$\underline{\modu}$ as triangulated categories, namely,
we have the following proposition.
\begin{prop}[\cite{KeVo, Ric1}]
  Let $A$ be a self-injective algebra. Then the following natural functor is an equivalence of triangulated categories,
  \begin{equation*}
    \mbox{$A$-$\underline{\modu}$}\rightarrow \DD_{\sg}(A).
  \end{equation*}
\end{prop}
Recall that for any $m\in\Z$, we define $$\HH_{\sg}^m(A, A):=\Hom_{\DD_{\sg}(A^e)}(A, A[m]).$$
Thanks to Corollary 6.4.1 in \cite{Bu}, we have the following descriptions for $\HH_{\sg}^*(A, A)$ in the case of a self-injective algebra $A$.
\begin{prop}[Corollary 6.4.1, \cite{Bu}]\label{prop-hom}
  Let $A$ be a self-injective algebra over a field $k$, denote $A^\vee:=\Hom_{A^e}(A, A^e)$.
  Then
  \begin{enumerate}
    \item $\HH_{\sg}^i(A, A)\cong \HH^i(A,A)$ for all $i>0$,
    \item $\HH_{\sg}^{-i}(A, A)\cong \Tor_{i-1}^{A^e}(A, A^{\vee})$ for all $i\geq 2$,
    \item there is an exact sequence
    \begin{equation*}
      \xymatrix{
      0\ar[r] & \HH^{-1}_{\sg}(A, A)\ar[r] & A^{\vee}\otimes_{A^e}A\ar[r] &  \Hom_{A^e}(A, A) \ar[r] & \HH^0_{\sg}(A, A)\rightarrow 0.
      }
    \end{equation*}
  \end{enumerate}
\end{prop}

\begin{rem}
Since $A$ is self-injective,  so is $A^e$. Hence we have for $i\geq 2$ and $n\geq 1$,
$$\Tor_{i-1}^{A^e}(A, A^{\vee})\cong \HH_{\sg}^{-i}(A, A)\cong\underline{\Hom}_{A^e}(\Omega^{i}(A), \Omega^{n+i}(A)\cong \Ext^{n}_{A^e}(A, \Omega^{n+i}(A)).$$
To simplify the notation, we denote each of these isomorphisms by $\lambda_{i, n}$.
In fact, we can write  the isomorphism $\lambda_{i, n}$ explicitly,
for $$\alpha\otimes a_1\otimes \cdots \otimes a_{i-1}\in\Tor_{i-1}^{A^e}(A, A^{\vee}),$$
 we have $$\lambda_{i, n}(\alpha\otimes a_1\otimes \cdots \otimes a_{i-1})\in\Ext^{n}_{A^e}(A, \Omega^{n+i}(A))$$ which is defined as follows,
\begin{equation}\label{equ-formular}
\begin{split}
\lambda_{i, n}(\alpha\otimes a_1\otimes \cdots \otimes a_{i-1})(b_1\otimes \cdots \otimes b_n)=\sum_j d(x_j\otimes a_{1, i-1}\otimes y_j\otimes b_{1, n}\otimes 1),
\end{split}
\end{equation}
where we write $\alpha(1):=\sum_jx_j\otimes y_j$. 
\end{rem}
Under the isomorphisms $\lambda_{i, n}$ above, we have the following relations between the  cap product and  cup product.
\begin{lemma}\label{lemma-cap-cup1}
Let $A$ be a self-injective algebra. Then we have the following commutative diagram
for $m\geq 1, n\geq 2$ and $n-m\geq 2$.
  \begin{equation*}
  \xymatrix{
  \HH^m(A, A)\otimes \HH^{-n}_{\sg}(A, A) \ar[r]^-{\cup}\ar[d]^{\id\otimes \lambda_n} & \HH_{\sg}^{m-n}(A, A)\ar[d]^{\lambda_{n-m}}\\
  \HH^m(A, A)\otimes \Tor^{A^e}_{n-1}(A, A^{\vee}) \ar[r]^-{\cap} & \Tor^{A^e}_{n-m-1}(A, A^{\vee}).
  }
\end{equation*}
\end{lemma}
\pf Take $f\otimes (\alpha\otimes a_1\otimes \cdots \otimes a_{n-1})\in \HH^m(A, A)\otimes \Tor^{A^e}_{n-1}(A, A^{\vee})$. Then
$$f\cap (\alpha\otimes a_1\otimes \cdots \otimes a_{n-1})=\sum_i x_i f(a_{1, m})\otimes y_i\otimes a_{m+1, n-1}\in\Tor^{A^e}_{n-m-1}(A, A^{\vee}),$$
where $\alpha(1):=\sum_i x_i\otimes y_i$.
Hence we have
$$\lambda_{n-m, m+1}(f\cap (\alpha\otimes a_1\otimes \cdots\otimes  a_{n-1})) \in \Ext^{m+1}(A, \Omega^{n+1}(A)).$$
From the formula in (\ref{equ-formular}), it follows that for any $(b_1\otimes \cdots\otimes b_{m+1})\in A^{\otimes  (m+1)}$,
 $$\lambda_{n-m, m+1}(f\cap (\alpha\otimes a_1\otimes \cdots \otimes a_{n-1}))(b_1\otimes \cdots\otimes  b_{m+1})=
\sum_id(x_if(a_{1, m})\otimes a_{m+1, n-1}\otimes y_i\otimes b_{1, m+1}\otimes 1).$$
We consider the following diagram,
\begin{equation}\label{equ-diagram1}
  \xymatrix{
 \Ext_{A^e}^m(A, A)\otimes \Ext_{A^e}^1(A, \Omega^{n+1}(A)) \ar[r]^-{\cup} &\Ext^{m+1}_{A^e}(A, \Omega^{n+1}(A))\\
  \HH^m(A, A)\otimes \HH^{-n}_{\sg}(A, A)\ar[u]_-{\lambda_{-m}\otimes \lambda_{n, 1}} \ar[r]^-{\cup}\ar[d]^{\id\otimes \lambda_n} & \HH_{\sg}^{m-n}(A, A)\ar[d]^{\lambda_{n-m}}\ar[u]_-{\lambda_{n-m, m+1}}\\
  \HH^m(A, A)\otimes \Tor^{A^e}_{n-1}(A, A^{\vee}) \ar@/^8pc/[uu]^{\id\otimes \lambda_{n, 1}} \ar[r]^-{\cap} & \Tor^{A^e}_{n-m-1}(A, A^{\vee})\ar@/_8pc/[uu]_{\lambda_{n-m, m+1}}
  }
\end{equation}
where $\cup$ in the first row represents the Yoneda product in the bounded derived category $\DD^b(A\otimes A^{\op})$.  It is clear that the top square in (\ref{equ-diagram1}) is commutative. Let us prove the commutativity of the outer square. For $$f\otimes (\alpha\otimes a_1\otimes \cdots a_{n-1})\in \HH^m(A, A)\otimes \Tor^{A^e}_{n-1}(A, A^{\vee}),$$
we have (via the up-right direction in Diagram (\ref{equ-diagram1}))
$$(\cup \circ (\id\otimes \lambda_{n, 1}))(f\otimes \alpha\otimes a_{1, n-1})\in \Ext_{A^e}^{m+1}(A, \Omega^{n+1}(A))$$
which sends $(b_{1}\otimes \cdots\otimes b_{m+1})\in A^{\otimes (m+1)}$ to
\begin{equation*}
(\cup \circ (\id\otimes \lambda_{n, 1}))(f\otimes \alpha\otimes a_{1, n-1})(b_{1, m+1})=\sum_i f(b_{1,m})d(x_i\otimes
a_{1, n-1}\otimes y_i\otimes b_{m+1}\otimes 1).
\end{equation*}
On the other hand, we have (via the right-up direction in Diagram (\ref{equ-diagram1}))
$$(\lambda_{n-m, m+1}\circ \cap)(f\otimes \alpha\otimes a_{1, n-1})\in  \Ext_{A^e}^{m+1}(A, \Omega^{n+1}(A))$$
which sends $(b_{1}\otimes \cdots\otimes b_{m+1})\in A^{\otimes (m+1)}$ to
\begin{equation*}
(\lambda_{n-m, m+1}\circ \cap)(f\otimes \alpha\otimes a_{1, n-1})(b_{1, m+1})=\sum_i d(x_if(a_{1, m})\otimes a_{m+1,
n-1}\otimes y_i\otimes b_{1, m+1}\otimes 1).
\end{equation*}
For $j=1, 2, \cdots, n,$ define
$$H_j\in \Hom_k(A^{\otimes m}, \Omega^{n+1}(A))$$ as follows,
\begin{equation*}
H_j(b_{1,m}):=\sum_i d (\id_{j}\otimes f)(x_i\otimes a_{1, n-1}\otimes y_i\otimes b_{1, m}\otimes 1).
\end{equation*}
By calculation, we have
\begin{equation*}
(\cup \circ (\id\otimes \lambda_{n, 1})-\lambda_{n-m, m+1}\circ \cap)(f\otimes \alpha\otimes a_{1, n-1})=
\sum_{j=1}^n (-1)^{\epsilon_j}\delta(H_j),
\end{equation*}
where $\epsilon_j\in \Z$ depends on $H_j$.
Hence in $\Ext_{A^e}^{m+1}(A, \Omega^{n+1}(A))$, we have
$$(\cup \circ (\id\otimes \lambda_{n, 1}))(f\otimes \alpha\otimes a_{1, n-1})=(\lambda_{n-m, m+1}\circ \cap)(f\otimes \alpha\otimes a_{1, n-1}).$$
So we have verified that the outer square in Diagram (\ref{equ-diagram1}) commutes, hence
the lower square also commutes.
 \epf

\begin{rem}
For the case $n-m=1$,
we have the following commutative diagram
\begin{equation*}
 \xymatrix{
  \HH^m(A, A)\otimes \HH^{-n}_{\sg}(A, A) \ar[r]^-{\cup}\ar[d]^{\id\otimes \lambda_n} & \HH_{\sg}^{1}(A, A)\ar@{^{(}->}[d]^{\lambda_{1}}\\
  \HH^m(A, A)\otimes \Tor^{A^e}_{n-1}(A, A^{\vee}) \ar[r]^-{\cap} & \Tor^{A^e}_{0}(A, A^{\vee}),
  }\end{equation*}
  where  the injection $\lambda_1:\HH_{\sg}^1(A, A)\rightarrow \Tor_0^{A^e}(A, A^{\vee})$ is defined in Proposition \ref{prop-hom}. Indeed, it is sufficient to prove that
  \begin{equation}\label{*}
    f\cap (\alpha\otimes a_{1, n-1})=\sum_ix_if(a_{1, m})\otimes y_i\in \Imm(\lambda_1)
  \end{equation} for
  $f\in \HH^m(A, A)$ and $\alpha\otimes a_{1, n-1}\in \Tor_{n-1}^{A^e}(A, A^{\vee}).$ From the exact sequence
  in Propostion \ref{prop-hom},
  \begin{equation*}
   \xymatrix{
      0\ar[r]& \HH^{-1}_{\sg}(A, A)\ar[r]^-{\lambda_1}  & A^{\vee}\otimes_{A^e}A\ar[r]^-{\mu_1} &  \Hom_{A^e}(A, A) \ar[r] & \HH^0_{\sg}(A, A)\rightarrow 0,
      }
        \end{equation*}
  it follows that (\ref{*}) is equivalent to
   \begin{equation*}
   \mu_1(f\cap (\alpha\otimes a_{1, n-1}))=\sum_ix_if(a_{1,m})y_i=0.
   \end{equation*}
   Indeed, we have
   \begin{equation*}
   \begin{split}
   \sum_ix_if(a_{1, m})y_i=&\sum_i f(x_ia_1\otimes a_{2, m})y_i+\sum_i \sum_{j=1}^{m-1}(-1)^{j}f(x_i\otimes a_{1, j-1}
   \otimes a_ja_{j+1}\otimes a_{j+2, m})y_i+\\
   &(-1)^mf(x_i\otimes a_{1, m-1})a_my_i\\
   =&0.
   \end{split}
   \end{equation*}
   Here we used $\delta(f)=0$ in the first identity and $d(\alpha\otimes a_{1, n-1})=0$ in the second identity.
\end{rem}

We are also interested in the case $n-m<0$. Before discussing this case, let us define some operators on Hochschild cohomology and homology. Let $A$ be an associative algebra (not necessarily, self-injective) over a field $k$.
We will define a (generalized) cap product as follows. Let $0\leq n<m$, define
\begin{equation*}
 \cap: C^m(A, A)\otimes C_{n}(A, A^{\vee})\rightarrow C^{m-n-1}(A, A)
\end{equation*}
$(f\cap (\alpha \otimes a_1\otimes \cdots\otimes a_n))(b_1\otimes \cdots \otimes b_{m-n-1}):=\sum_ix_if(a_1\otimes \cdots \otimes a_n\otimes y_i\otimes b_1\otimes \cdots \otimes b_{m-n-1}),$
where $\alpha(1):=\sum_ix_i\otimes y_i$.
\begin{lemma}\label{lemma-cap}
The $\cap$ induces a well-defined operator on the level of homology,
$$\cap: \HH^m(A, A)\otimes \HH_n(A, A^{\vee})\rightarrow \HH^{m-n-1}(A, A).$$
\end{lemma}
\pf It is sufficient to check that
\begin{equation*}
\begin{split}
Z^m(A, A)\cap Z_n(A, A^{\vee})\subset Z^{m-n-1}(A, A),\\
Z^m(A, A)\cap B_n(A, A^{\vee})\subset B^{m-n-1}(A, A),\\
B^m(A, A)\cap Z_n(A, A^{\vee})\subset B^{m-n-1}(A, A).\\
\end{split}
\end{equation*}
Let $f\in Z^m(A, A)$ and $z=\sum \alpha\otimes a_{ 1}\otimes \cdots \otimes a_{n}\in Z_n(A, A^{\vee}).$
Then we have
\begin{equation*}
\begin{split}
&\delta(f\cap z)(b_1\otimes \cdots \otimes b_{m-n})\\
=&b_1(f\cap z)(b_{2, m-n})+\sum_{j=1}^{m-n-1}(-1)^j(f\cap z)(b_{1, j-1}\otimes b_jb_{j+1}\otimes b_{j+2, m-n})+\\
&(-1)^{m-n}(f\cap z)(b_{1, m-n-1})b_{m-n}\\
=&b_1\sum_k x_{ k}f(a_{1, n}\otimes y_k\otimes b_{2, m-n})+\sum_{j=1}^{m-n-1}(-1)^j\sum_k x_{ k}f(a_{1, n}\otimes y_k\otimes b_{2, j-1}\otimes b_jb_{j+1}\otimes b_{j+2, m-n})+\\
&(-1)^{m-n}\sum_k x_{ k}f(a_{1, n}\otimes y_k\otimes b_{1, m-n-1})b_{m-n}\\
=&0,
\end{split}
\end{equation*}
where the last identity follows from the fact that $\delta(f)=0$ and $dz=0$.
Similarly, for $f\in\Hom(A^{\otimes m-1}, A)$ and $z\in Z_n(A, A^{\vee})$
we have $$\delta(f)\cap z=\delta(f\cap z),$$
and for $f\in Z^m(A, A)$ and $z\in C_{n+1}(A, A^{\vee})$
$$f\cap \partial_{n+1}(z)=\delta(f\cap z).$$
\epf

Let $A$ be an associative algebra.
We will also define a (generalized) cup product on $\Tor_*^{A^e}(A, A^{\vee})$. Let $m, n\in \Z_{\geq 0},$
\begin{equation}\label{gene-cup}
\cup:C_{m}(A, A^{\vee})\times C_{n}(A, A^{\vee})\rightarrow C_{m+n+1}(A, A^{\vee})
\end{equation}
is defined as follows,
take $$\sum_i (x_i\otimes y_i)\otimes a_{1, m}\in C_{m}(A, A^{\vee})$$ and
$$\sum _j (x_j'\otimes y_j')\otimes b_{1,n}\in C_{n}(A, A^{\vee}),$$
we define
\begin{equation*}
(\sum_i x_i\otimes y_i\otimes a_{1, m})\cup(\sum_j x_j'\otimes y_j'\otimes b_{1,n}):=\sum_{i, j} (x_i\otimes y_j'y_i)\otimes a_{1,m}\otimes x_j'\otimes b_{1,n}.
\end{equation*}
\begin{lemma}
The generalized cup product defined above is well-defined on $\Tor_{> 0}^{A^e}(A, A^{\vee})$. Moreover, it is graded commutative.
\end{lemma}
\pf Let $$\alpha=\sum (x_i\otimes y_i)\otimes a_{1, m}\in Z_{m}(A, A^{\vee})$$ and $$\beta=\sum (x_j'\otimes y_j')\otimes b_{1,n}\in Z_{n}(A, A^{\vee}),$$
then we have the following,
\begin{equation*}
\begin{split}
d(\alpha\cup\beta)=&d(\sum (x_i\otimes y_j'y_i)\otimes a_{1,m}\otimes x_j'\otimes b_{1,n})\\
=&\sum (x_ia_1\otimes y_j'y_i)\otimes a_{2,m}\otimes x_j'\otimes b_{1,n})+\\
&\sum^{m-1}_{k=1}(-1)^{k} \sum (x_i\otimes y_j'y_i)\otimes a_{1,i-1}\otimes a_ia_{i+1}\otimes a_{i+2, m}\otimes x_j'\otimes b_{1,n}+\\
&\sum (-1)^{m} \sum (x_i\otimes y_j'y_i)\otimes a_{1, m-1}\otimes a_mx_j'\otimes b_{1,n}+\\
&\sum(-1)^{m+1}\sum (x_i\otimes y_j'y_i)\otimes a_{1, m}\otimes x_j'b_1\otimes b_{2,n}+\\
&\sum\sum^{n-1}_{k=1}(-1)^{m+1+k}(x_i\otimes y_j'y_i)\otimes a_{1, m}\otimes x_j'\otimes b_{1,k-1}\otimes b_kb_{k+1}\otimes b_{k+2,n}+\\
&\sum (-1)^{m+n+1}(x_i\otimes b_ny_j'y_i)\otimes a_{1, m}\otimes x_j'\otimes b_{1,n-1}\\
=&0,
\end{split}
\end{equation*}
where the last identity follows from the fact that $d\alpha=d\beta=0$.
Hence we have
\begin{equation*}
Z_m(A, A^{\vee})\cup Z_{n}(A, A^{\vee})\subset Z_{m+n+1}(A, A^{\vee}).
\end{equation*}
Similarly, the followings can be verified
\begin{equation*}
\begin{split}
Z_m(A, A^{\vee})\cup B_{n}(A, A^{\vee})\subset B_{m+n+1}(A, A^{\vee}),\\
B_m(A, A^{\vee})\cup Z_{n}(A, A^{\vee})\subset B_{m+n+1}(A, A^{\vee}).
\end{split}
\end{equation*}
It remains to verify the graded commutativity. By calculation, we have
\begin{equation*}
\begin{split}
\alpha\cup \beta-(-1)^{mn}\beta\cup \alpha=\sum \sum_{k=0}^{m}(-1)^{m+n(k-1)}d((x_i\otimes y_i)
\otimes a_{1, m-k}\otimes x_j'\otimes b_{1,n}\otimes y_j'\otimes a_{m-k+1,m}).
\end{split}
\end{equation*}
Hence, we have $$\alpha\cup \beta-(-1)^{mn}\beta\cup \alpha\in B_{m+n+1}(A, A^{\vee}).$$
Thus, $\alpha\cup\beta=(-1)^{mn}\beta\cup \alpha$. Therefore, we have finished the proof.
\epf

Let us go back to our special case, then we have the following.
\begin{lemma}\label{lemma-cap-cup2}
Let $A$ be a self-injecitve algebra over a field $k$. Then the following diagram commutes for
$m\geq 1, n\geq 2$ and $m-n\geq 1,$
\begin{equation*}
\xymatrix{
\HH^m(A, A)\otimes \HH_{\sg}^{-n}(A, A)\ar[r]^-{\cup} & \HH_{\sg}^{m-n}(A, A)\\
\HH^m(A, A)\otimes \Tor^{A^e}_{n-1}(A, A^{\vee})\ar[r]^-{\cap} \ar[u]^{\id\otimes \lambda_{n}}&\HH^{m-n}(A, A)\ar[u]^{\cong}.
}
\end{equation*}
\end{lemma}
\pf The proof is similar to the proof of Lemma \ref{lemma-cap-cup1}. \epf
\begin{rem}
For $m-n=0,$ we have the following commutative diagram,
\begin{equation*}
\xymatrix{
\HH^m(A, A)\otimes \HH_{\sg}^{-m}(A, A)\ar[r]^-{\cup} & \HH_{\sg}^{0}(A, A)\\
\HH^m(A, A)\otimes \Tor^{A^e}_{m-1}(A, A^{\vee})\ar[r]^-{\cap} \ar[u]^{\id\otimes \lambda_{n}}&\HH^{0}(A, A)\ar@{->>}[u]^{\pi_0},
}
\end{equation*}
where $\pi_0: \HH^0(A, A)\rightarrow \HH_{\sg}^0(A, A)$ is the surjection defined in Proposition \ref{prop-hom}.
\end{rem}

Similarly, we also have the following lemma.
\begin{lemma}\label{lemma-negative-cup}
Let $A$ be a self-injective algebra over a field $k$. Then we have the following commutative diagram for $m\geq 2,
n\geq 2$,
\begin{equation*}
\xymatrix{
\HH_{\sg}^{-m}(A, A)\otimes \HH_{\sg}^{-n}(A, A)\ar[r]^-{\cup} & \HH_{\sg}^{-m-n}(A, A)\\
\Tor_{m-1}^{A^e}(A, A^{\vee})\otimes \Tor^{A^e}_{n-1}(A, A^{\vee})\ar[r]^-{\cup} \ar[u]^{\lambda_m\otimes \lambda_{n}}&\Tor_{m+n-1}^{A^e}(A, A^{\vee})\ar[u]_{\lambda_{m+n}}.
}
\end{equation*}
\end{lemma}
\begin{rem}\label{rem-negative-cup}
For $m=n=1$,  we have the following commutative diagram,
\begin{equation*}
\xymatrix{
\HH_{\sg}^{-1}(A, A)\otimes \HH_{\sg}^{-1}(A, A)\ar[r]^-{\cup} \ar[d]_-{(\lambda_1\otimes \lambda_1)^{-1}}& \HH_{\sg}^{-2}(A, A)\\
\Tor_{0}^{A^e}(A, A^{\vee})\otimes \Tor^{A^e}_{0}(A, A^{\vee})\ar[r]^-{\cup} &\Tor_{1}^{A^e}(A, A^{\vee})\ar[u]_{\lambda_{m+n}}.
}
\end{equation*}
\end{rem}

Therefore, in conclusion,  the graded commutative associative algebra  structure on
$(\HH^*_{\sg}(A, A), \cup)$ becomes well-understood. Namely, it can be interpreted as the (generalized) cap product and the (generalized) cup product. Next we will investigate the graded Lie algebra structure on  $(\HH^*_{\sg}(A, A), [\cdot,\cdot])$ in the case of a symmetric algebra $A$.

Before starting the case of symmetric algebras, let us recall the Connes B-operator on Hochschild homology and the structure of a Batalin-Vilkovisky (BV) algebra. For more details, we refer to \cite{Con, Lod,Xu}.

\begin{defn}
Let $A$ be an associative algebra over a commutative algebra $k$. We define an operator on Hochschild homology,
$$B: C_r(A, A)\rightarrow C_{r+1}(A, A),$$
which sends $a_0\otimes \cdots \otimes a_r\in C_r(A, A)$
to
\begin{equation*}
\begin{split}
B(a_0\otimes \cdots \otimes a_r):&=\sum_{i=0}^r(-1)^{ir}1\otimes a_i\otimes \cdots \otimes a_r\otimes a_0\otimes \cdots \otimes a_{i-1}+\\
&(-1)^{ir}a_i\otimes 1\otimes a_{i+1}\otimes \cdots \otimes a_r\otimes a_0\otimes \cdots a_{i-1}.
\end{split}
\end{equation*}
It is easy to check that $B$ is a chain map satisfying $$B\circ B=0,$$
which induces an operator (still denote by B),
$$ B : \HH_r(A, A)\rightarrow \HH_{r+1}(A, A).$$
We call this operator  Connes B-operator.
\end{defn}
\begin{rem}
We also consider the Connes B-operator on normalized Hochschild complex.
$$\overline{B}: \overline{C}_r(A, A)\rightarrow \overline{C}_{r+1}(A, A),$$ which sends
$a_{0, r}\in\overline{C}_r(A, A)$ to
$$\overline{B}(a_{0, r}):=\sum_{i=0}^r(-1)^{ir}1\otimes a_i\otimes \cdots \otimes a_r\otimes a_0\otimes \cdots \otimes a_{i-1}.$$
\end{rem}

\begin{defn}\label{defn-BV}
A Batalin-Vilkovisky algebra (BV algebra for short) is a Gerstenhaber algebra $(\calH^*, \cup, [\cdot, \cdot])$ together with an operator
$\Delta: \calH^*\rightarrow \calH^{*-1}$ of degree $-1$ such that $\Delta\circ \Delta =0, \Delta(1)=0$ and satisfying the following BV identity,
$$[\alpha, \beta]=(-1)^{|\alpha|+1}\Delta(\alpha\cup \beta)+(-1)^{|\alpha|}\Delta(\alpha)\cup \beta+\alpha\cup \Delta(\beta)$$
for homogeneous elements $\alpha, \beta \in\calH^*$.
\end{defn}


From here onwards, assume that $k$ is a field. Let  $A$ be a symmetric $k$-algebra (i.e. there is a symmetric, associative and non-degenerate inner product
$\langle\cdot, \cdot \rangle: A\otimes A\rightarrow k$).
Then the inner product $\langle\cdot, \cdot\rangle$ induces an $A$-$A$-bimodule
isomorphism
\begin{equation}\label{equ-t}
 \begin{tabular}{rccc}
   $t: $ & $A$ & $\rightarrow$ &$ D(A):=\Hom_k(A, k)$\\
    & $a $ & $\mapsto$ & $\langle a, -\rangle.$
    \end{tabular}
\end{equation}
where the $A$-$A$-bimodule structure on $D(A)$ is given as follows,
for $f\in D(A)$ and $a\otimes b\in A\otimes A^{\op}$, $$((a\otimes b)f)(c)=f(cba).$$
This isomorphism $t$ induces the following isomorphism
\begin{equation*}
  \begin{tabular}{rccc}
    $t\otimes \id: $ & $A\otimes A$ & $\rightarrow$ &$ D(A)\otimes A \cong\End(A)$\\
    & $a\otimes b $ & $\mapsto$ & $t(a)\otimes b\mapsto ( x\mapsto t(a)(x) b).$
    \end{tabular}
\end{equation*}
We define the element $$(t\otimes \id)^{-1}(\id):=\sum_i e_i\otimes f_i\in A\otimes A$$ as the Casimir element of $A$ (with respect to the inner product $\langle \cdot, \cdot \rangle$) (cf. \cite{Brou}).
The following proposition states some properties on Casimir element.
\begin{prop}[Proposition 3.3. \cite{Brou}]\label{prop-broue}
  \begin{enumerate}
    \item For all $a, a'\in A$, we have
    $$\sum_i ae_ia'\otimes f_i=e_i\otimes a'f_ia.$$
    \item The map
    \begin{equation*}
      \begin{tabular}{rccc}
     $A$ & $\rightarrow$ &$ \Hom_{A^e}(A, A^e)$\\
     $a $ & $\mapsto$ & $\sum_i e_ia\otimes f_i.$
    \end{tabular}
    \end{equation*}
    is an right $A$-$A$-bimodule isomorphism. Here $A$ is a right $A$-$A$-bimodule defined as follows, for $a\in A$ and $b\otimes c\in A^e$, $$a\cdot(b\otimes c):=cab,$$
    the right $A$-$A$-bimodule structure on $\Hom_{A^e}(A, A^e)$ is given by,
    for $f\in \Hom_{A^e}(A, A^e)$ and $b\otimes c\in A^e$,
    $$f\cdot(b\otimes c)(a):=f(cba),$$
    and we identify $\Hom_{A^e}(A, A^e)$ as $$(A\otimes A)^A:=\{\sum a_i\otimes b_i\in A\otimes A \ | \  \sum aa_i\otimes b_i=\sum a_i\otimes b_ia, \ \mbox{for any $a\in A$} \}.$$
  \end{enumerate}
\end{prop}

Since we have the following isomorphisms via the isomorphism $t$ defined in (\ref{equ-t}) above,
\begin{equation*}
  \begin{split}
    D(C_n(A, A))&\cong D(A\otimes_{A^e}D^2(A^{\otimes n+2}))\cong \Hom_{A^e}(A, D(A^{\otimes n+2}))\\
    &\cong \Hom_{A^e}(D^2(A^{\otimes n+2}), D(A))\cong  C^n(A, A)
      \end{split}
\end{equation*}
where the third isomorphism follows from the fact that $D$ induces an equivalence between 
$A^e$-$\modu$ and $(A^e$-$\modu)^{\op}$ and the forth isomorphism is induced from the isomorphism
$t$.
Hence we have a duality between Hochschild homology and cohomology, for $n\in Z_{\geq 0}$,
\begin{equation}\label{equ-dua}
\HH_n(A, A)^*\cong \HH^n(A, A).
\end{equation}
Hence from Proposition \ref{prop-hom} and \ref{prop-broue},
 we have for $n\geq 1$,
 $$\kappa_{n}:\Ext^{1}(A, \Omega^{n+2}(A))\cong \HH_{\sg}^{-n-1}(A, A)\cong\Tor_n^{A^e}(A, A^{\vee})\cong \Tor_n^{A^e}(A, A))\cong \HH^n(A, A)^*.$$
 In the rest of our paper, for simplicity, we often use the same $\kappa_{n}$ to indicate any of those natural isomorphisms above. For example, we have the following isomorphism,
 \begin{equation*}
   \begin{tabular}{ccccc}
    $ \kappa_{n}:$& $\Tor_{n}^{A^e}(A, A)$ & $\rightarrow$ & $\Ext^{1}(A, \Omega^{n+2}(A))$\\
    & $a_0\otimes a_1\otimes \cdots \otimes a_n$ & $\mapsto$ & $(b\mapsto \sum_i d (e_ia_0\otimes a_{1, n}\otimes f_i\otimes b\otimes 1)).$
   \end{tabular}
 \end{equation*}
Moreover we have the following result on symmetric algebras.
\begin{thm}[\cite{Tra}, \cite{Men}]\label{thm-tra}
Let $A$ be a symmetric algebra over a field $k$. Then $$(\HH^*(A, A), \cup, [\cdot, \cdot], \Delta)$$ is a BV algebra, where the BV-operator $\Delta$ is the dual of the Connes $ B $-operator via the duality (\ref{equ-dua}) above.
\end{thm}

Now let us state our propositions.
\begin{prop}\label{prop-tor1}
Let $A$ be a symmetric algebra over a field $k$. Then we have the following commutative diagram for $m\geq 1, n\geq 2$ and $n-m\geq 1,$
\begin{equation*}
\xymatrix{
\HH^m(A, A)\otimes \Ext^1_{A^e}(A, \Omega^{n+1}(A)) \ar[r]^-{[\cdot, \cdot]} & \Ext^m_{A^e}(A, \Omega^{n+1}(A))\\
\HH^m(A, A)\otimes \HH_{\sg}^{-n}(A, A)\ar[r]^-{[\cdot, \cdot]} \ar[u]^{\cong} & \HH_{\sg}^{m-n-1}(A, A)\ar[u]_{\cong}\\
\HH^m(A, A)\otimes \Tor_{n-1}^{A^e}(A, A)\ar@/^8pc/[uu]^{\id\otimes \kappa_n}\ar[r]^-{\{\cdot, \cdot\}} \ar[u]^{\id\otimes \kappa_n} & \Tor_{n-m}^{A^e}(A, A)\ar[u]_{\kappa_{n-m+1}}\ar@/_8pc/[uu]_-{\kappa_{n-m+1}}
}
\end{equation*}
where $\{\cdot, \cdot\}$ is defined as follows, for $f\in\HH^m(A, A), \alpha\in\Tor_{n-1}^{A^e}(A, A)$,
\begin{equation*}
\begin{split}
\{f, \alpha\}:=(-1)^{m}\Delta(f)\cap \alpha+f\cap  B (\alpha)+(-1)^{m+1} B (f\cap\alpha).
\end{split}
\end{equation*}
\end{prop}
\pf  From the definition of the Gerstenhaber bracket $[\cdot,\cdot]$, it follows that the top square is commutative. Hence it remains to check
the commutativity of the outer square.
Let $f\in \HH^m(A, A)$ and $$z:=\sum a_0\otimes a_1\otimes \cdots \otimes a_{n-1}\in \Tor_{n-1}^{A^e}(A, A).$$ Then
\begin{equation*}
  \begin{split}
   & \kappa_{n-m+1}(\{f, z\})(b_1\otimes \cdots \otimes b_m)\\
    =&\sum (-1)^{m+(m-1)(n-1)} d(e_ja_0\Delta(f)(a_{1, m-1})\otimes a_{m, n-1}\otimes f_j\otimes b_{1, m}\otimes 1)+\\
     &\sum \sum_{i=0}^{n-1}(-1)^{i(n-1)+mn}d(e_j(f\otimes \id)(a_{i, n-1}\otimes a_{0, i-1})\otimes f_j\otimes  b_{1, m}\otimes 1)+\\
     &\sum (-1)^{m+1+m(n-1)}d(e_j\otimes a_0f(a_{1, m})\otimes a_{m+1, n-1}\otimes f_j\otimes b_{1, m}\otimes 1)+\\
     &\sum \sum_{i=1}^{n-m-1}(-1)^{m+1+m(n-1)+i(n-m-1)}d(e_j\otimes a_{i+m, n-1}\otimes a_0f(a_{1, m})\otimes a_{m+1, i+m-1}\otimes f_j\otimes b_{1, m}\otimes 1)\\
  \end{split}
\end{equation*}
\begin{equation*}
  \begin{split}
   & [f, \kappa_n(z)](b_1\otimes \cdots \otimes b_m)=f\bullet \kappa_n(z)(b_{1, m})-\kappa_n(z)\circ f(b_{1, m})\\
   =&\sum \sum_{i=1}^m(-1)^{(m+1-i)(n+1)}d(f\otimes \id)(b_{1, i-1}\otimes d(e_ja_0\otimes a_{1, n-1}\otimes f_j\otimes b_i\otimes 1)\otimes b_{i+1, m})\otimes 1+\\
   &\sum_{i=1}^{n+1}(-1)^{n+1+(n+1-i)m}d(id_i\otimes f)(d(e_ja_0\otimes a_{1, n-1}\otimes f_j\otimes b_1\otimes 1)\otimes b_{2, m}\otimes 1)-\\
   & \sum d(e_ja_0\otimes a_{1, n-1}\otimes f_j\otimes f(b_{1, m})\otimes 1)\\
   =&\sum \sum_{i=1}^m(-1)^{(m-i)(n+1)}d(f\otimes \id)(b_{1, i-1}\otimes e_ja_0\otimes a_{1, n-1}\otimes f_j\otimes b_{i, m})\otimes 1+\\
   &\sum \sum_{i=1}^{n}(-1)^{(n+1-i)m}d(id_i\otimes f)(e_ja_0\otimes a_{1, n-1}\otimes f_j\otimes b_{1, m}\otimes 1)
  \end{split}
\end{equation*}
Claim that
\begin{equation*}
\begin{split}
0=&\sum \sum_{i=1}^{n-m+1}(-1)^{i(n-1)+mn}d(e_j(f\otimes \id)(a_{i, n-1}\otimes a_{0, i-1})\otimes f_j\otimes  b_{1, m}\otimes 1)+\\
&\sum (-1)^{m+1+m(n-1)}d(e_j\otimes a_0f(a_{1, m})\otimes a_{m+1, n-1}\otimes f_j\otimes b_{1, m}\otimes 1)+\\
     &\sum \sum_{i=1}^{n-m-1}(-1)^{m+1+m(n-1)+i(n-m-1)}d(e_j\otimes a_{i+m, n-1}\otimes a_0f(a_{1, m})\otimes a_{m+1, i+m-1}\otimes f_j\otimes b_{1, m}\otimes 1)-\\
&\sum \sum_{i=1}^{n-m}(-1)^{(n+1-i)m}d(id_i\otimes f)(e_ja_0\otimes a_{1, n-1}\otimes f_j\otimes b_{1, m}\otimes 1)
\end{split}
\end{equation*}
Hence we have
\begin{equation}\label{equ-tor1}
\begin{split}
&( \kappa_{n-m+1}(\{f, z\})-[f, \kappa_n(z)])(b_1\otimes \cdots \otimes b_m)\\
=&\sum (-1)^{m+(m-1)(n-1)} d(e_ja_0\Delta(f)(a_{1, m-1})\otimes a_{m, n-1}\otimes f_j\otimes b_{1, m}\otimes 1)+\\
&\sum(-1)^{mn}d(e_jf(a_{0, m-1})\otimes a_{m, n-1}\otimes f_j\otimes b_{1, m}\otimes 1+\\
     &\sum \sum_{i=n-m+2}^{n-1}(-1)^{i(n-1)+mn}d(e_j(f\otimes \id)(a_{i, n-1}\otimes a_{0, i-1})\otimes f_j\otimes  b_{1, m}\otimes 1)-\\
     &\sum \sum_{i=1}^m(-1)^{(m-i)(n+1)}d(f\otimes \id)(b_{1, i-1}\otimes e_ja_0\otimes a_{1, n-1}\otimes f_j\otimes b_{i, m})\otimes 1 -\\
     & \sum \sum_{i=n-m+1}^{n}(-1)^{(n+1-i)m}d(id_i\otimes f)(e_ja_0\otimes a_{1, n-1}\otimes f_j\otimes b_{1, m}\otimes 1)  
     \end{split}
\end{equation}

We have the following
\begin{equation}\label{equ-tor2}
  \begin{split}
  &\sum_j d(\id_n\otimes f)(e_ja_0\otimes a_{1, n-1}\otimes  f_j\otimes b_{1, m}\otimes 1)\\
    =& \sum_{j, k}d(e_ja_0\otimes a_{1, n-1} \otimes \langle e'_kf(f_j\otimes b_{1, m-1}), 1\rangle f'_k\otimes b_m\otimes 1)\\
 =  &\sum_{j, k}\sum_{l=1}^{m-1}\sum_{i=l}^{m-1}(-1)^{(m+l)(m-i)+1}\delta(d(e_ja_0\otimes a_{1,n-l}\otimes  \langle f(\id_{m-i-1}\otimes e_k'\otimes a_{n-l+1, n-1}\otimes f_j\otimes \id_{i-l}), 1\rangle f_k'\\
 &\otimes \id_{l}\otimes 1)(b_{1, m})+\sum_{j}\sum_{i=1}^{m}(-1)^{m(m-1)}d(e_ja_0\otimes a_{1, n-1}\otimes \Delta(f)(b_{1, m-1})f_j\otimes b_m\otimes 1)+\\
 &\sum_{k}\sum_{l=1}^m(-1)^{(m-1)(m+l)}d(f(b_{1, m-l}\otimes e_k'\otimes a_{n-l+1, n-1})a_0\otimes a_{1, n-l}\otimes f_k'\otimes b_{m-l+1, m}\otimes 1)+\\
 & \sum_{j}\sum_{l=1}^{m-1}(-1)^{ml+1}d(e_ja_0\otimes a_{1, n-l-1}\otimes f(a_{n-l, n-1}\otimes f_j\otimes b_{1, m-l-1})\otimes b_{m-l, m}\otimes 1\\
      \end{split}
\end{equation}
Combining (\ref{equ-tor1}) with (\ref{equ-tor2}), we obtain that
\begin{equation}\label{equ-tor3}
\begin{split}
&( \kappa_{n-m+1}(\{f, z\})-[f, \kappa_n(z)])(b_1\otimes \cdots \otimes b_m)\\
=&\sum(-1)^{mn}d(e_jf(a_{0, m-1})\otimes a_{m, n-1}\otimes f_j\otimes b_{1, m}\otimes 1+\\
     &\sum \sum_{i=n-m+1}^{n-1}(-1)^{i(n-1)+mn}d(e_j(f\otimes \id)(a_{i, n-1}\otimes a_{0, i-1})\otimes f_j\otimes  b_{1, m}\otimes 1)-\\
     &\sum \sum_{i=1}^m(-1)^{(m-i)(n+1)}d(f\otimes \id)(b_{1, i-1}\otimes e_ja_0\otimes a_{1, n-1}\otimes f_j\otimes b_{i, m})\otimes 1 -\\
     & \sum_{j}\sum_{l=1}^m(-1)^{(m-1)(m+l)+m}df(b_{1, m-l}\otimes e_j\otimes a_{n-l+1, n-1})a_0\otimes a_{1, n-l}\otimes f_j\otimes b_{m-l+1, m}\otimes 1
          \end{split}
\end{equation}
Let us compute the following term in (\ref{equ-tor3}).
\begin{equation}\label{equ-tor4}
\begin{split}
&\sum d(f\otimes \id)(b_{1, m-1}\otimes e_ja_0\otimes a_{1, n-1}\otimes f_j\otimes b_{m})\otimes 1\\
=&\sum df(b_{1,m-1}\otimes e_ja_0)\otimes a_{1, n-1}\otimes f_j\otimes b_m\otimes 1\\
=&(-1)^{m-1}\sum db_1f(b_{2, m-1}\otimes e_j\otimes a_0)\otimes a_{1, n-1}\otimes f_j\otimes b_m\otimes 1+\\
&\sum\sum_{i=1}^{m-2}(-1)^{m-1+i}f(b_{1, i-1}\otimes b_ib_{i+1}\otimes b_{i+2, m-1} \otimes e_j\otimes
a_0)\otimes a_{1, n-1}\otimes f_j\otimes b_m\otimes  1+\\
&f(b_{1, m-2}\otimes b_{m-1}e_j\otimes a_0)\otimes a_{1, n-1}\otimes f_j\otimes b_m\otimes 1+df(b_{1, m-1}\otimes e_j)a_0\otimes a_{1, n-1}\otimes f_j\otimes b_m\otimes 1\\
=&\sum_{i=1}^{m-1}(-1)^{m-1+(i-1)n}\sum\delta(df(\id_{m-i-1}\otimes e_j\otimes a_{0, i-1})\otimes a_{i, n-1}\otimes \id_i\otimes 1)(b_{1, m})\\
&\sum_{i=1}^{m-1}(-1)^{(n-1)i+1}df(b_{1, m-1-i}\otimes e_ja_0\otimes a_{1, i})\otimes a_{i+1, n-1}\otimes f_j\otimes b_{m-i, m}\otimes 1+\\
&\sum_{i=1}^{m}(-1)^{(m-1)i+m+1}df(b_{1, m-i}\otimes e_j\otimes a_{n-i+1, n-1})a_0\otimes a_{1, n-i}\otimes f_j\otimes b_{m-i+1, m}\otimes 1+\\
&\sum (-1)^{mn}de_jf(a_{0, m-1})\otimes a_{m, n-1}\otimes f_j\otimes b_{1, m}\otimes 1+\\
&\sum_{i=n-m+1}^{n-1}(-1)^{i(n-1)+mn}de_j(f\otimes \id)(a_{i, n-1}\otimes a_{0, i-1})\otimes f_j\otimes b_{1, m}\otimes 1
\end{split}
\end{equation}
From (\ref{equ-tor4}), it follows that
$$ \kappa_{n-m+1}(\{f, z\})-[f, \kappa_n(z)]\in Z^m(A, \Omega^{n+1}(A)),$$
hence we have the following identity in $\Ext^m_{A^e}(A, \Omega^{n+1}(A))$
$$ \kappa_{n-m+1}(\{f, z\})=[f, \kappa_n(z)]).$$
Therefore we have completed our proof.
\epf

\begin{rem}\label{rem-bracket}
We also have the following commutative diagram for $m\geq 1,  n\geq 2$ and $n-m\geq 1$,
\begin{equation*}
\xymatrix{
\HH^m(A, A)\otimes \HH^{n-1}(A, A)^*\ar[r]^-{[\cdot, \cdot]^*} \ar[d]^-{\id\otimes \kappa_{n}} & \HH^{n-m}(A, A)^*\ar[d]^{\kappa_{m-n+1}}\\
\HH^m(A, A)\otimes \Tor_{n-1}^{A^e}(A, A)\ar[r]^-{\{\cdot, \cdot\}} & \Tor_{n-m}^{A^e}(A, A)
}
\end{equation*}
where $[\cdot,\cdot]^*$ is defined as follows, for any $f\in \HH^m(A, A)$ and $\alpha\in \HH^{n-1}(A, A)^*$,
$$[f, \alpha]^*(-):=\langle\alpha, [f, -]\rangle.$$
In fact, for $f\in\HH^m(A, A), \alpha\in\Tor_{n-1}^{A^e}(A, A)$ and $g\in \HH^{n-m}(A, A)$,
\begin{equation*}
\begin{split}
[f, \kappa_{n}^{-1}(\alpha)]^*(g)=&\langle \kappa_{n}^{-1}(\alpha), [f, g]\rangle\\
=&\langle\kappa_{n}^{-1}(\alpha), (-1)^{m}\Delta(f)\cup g+f\cup\Delta(g)+(-1)^{m+1}\Delta(f\cup g)\rangle\\
=&(-1)^m\langle\kappa_{m-n+1}^{-1}(\Delta(f)\cap\alpha), g \rangle+\langle \kappa_{m-n+1}^{-1}(B (f\cap\alpha)), g \rangle+\\
&(-1)^{m+1}\langle \kappa_{m-n+1}^{-1}(f\cap B (\alpha)), g \rangle\\
=&\langle \kappa_{m-n+1}^{-1}(\{f, \alpha\} ), g \rangle
\end{split}
\end{equation*}
where the second identity is the BV identity (cf. Definition \ref{defn-BV}), hence it follows that the diagram above is commutative.
\end{rem}

Similarly, we obtain the following proposition.
\begin{prop}\label{prop-ext1}
Let $A$ be a symmetric algebra over a field $k$. Then we have the following commutative diagram for $m\geq 1, n\geq 2$ and $n-m\leq -2,$
\begin{equation*}
\xymatrix{
\HH^m(A, A)\otimes \Ext^{1}_{A^e}(A, \Omega^{n+1}(A))\ar[r] &\Ext^m_{A^e}(A, \Omega^{n+1}(A))\\
\HH^m(A, A)\otimes \HH_{\sg}^{-n}(A, A)\ar[u]\ar[r]^-{[\cdot, \cdot]}  & \HH_{\sg}^{m-n-1}(A, A)\ar[u]\\
\HH^m(A, A)\otimes \Tor_{n-1}^{A^e}(A, A)\ar[u]^{\id\otimes \kappa_n}\ar[r]^-{\{\cdot, \cdot\}}  & \HH^{m-n-1}(A, A)
\ar[u]_{\kappa_{n-m+1}}
}
\end{equation*}
where $\{\cdot, \cdot\}$ is defined as follows: for any $f\in\HH^m(A, A)$ and $\alpha\in \Tor_{n-1}^{A^e}(A, A)$,
\begin{equation*}
\{f, \alpha\}:=(-1)^{m}\Delta(f)\cap \alpha+f\cap  B (\alpha)+(-1)^{m+1} \Delta (f\cap\alpha),
\end{equation*}
where $\alpha$ represents the generalized cap product defined in Lemma \ref{lemma-cap}.
\end{prop}
\pf The proof is similar to the one of Proposition \ref{prop-tor1}. Take $f\in\HH^m(A, A)$ and
$$z:=\sum a_0\otimes a_{1, n-1}\in\Tor_{n-1}^{A^e}(A, A).$$ Then we have
\begin{equation}\label{ext-0}
\begin{split}
&\kappa_{n-m+1}(\{f, z\})(b_{1, m})\\
=&(-1)^m \sum de_ja_0\Delta(f)(a_{1, n-1}\otimes f_j\otimes b_{1, m-n-1})\otimes b_{m-n, m}\otimes 1+\\
&\sum \sum_{i=0}^{n-1}(-1)^{i(n-1)}de_jf(a_{i, n-1}\otimes a_{0, i-1}\otimes f_j\otimes b_{1, m-n-1})\otimes b_{m-n, m}\otimes 1+\\
&(-1)^{m+1}d\Delta(f\cap z)(b_{1, m-n-1})\otimes b_{m-n, m}\otimes 1\\
&\\
   & [f, \kappa_n(z)](b_1\otimes \cdots \otimes b_m)\\
     =&\sum \sum_{i=1}^m(-1)^{(m-i)(n+1)}d(f\otimes \id)(b_{1, i-1}\otimes e_ja_0\otimes a_{1, n-1}\otimes f_j\otimes b_{i, m})\otimes 1+\\
   &\sum \sum_{i=1}^{n}(-1)^{(n+1-i)m}d(\id_i\otimes f)(e_ja_0\otimes a_{1, n-1}\otimes f_j\otimes b_{1, m}\otimes 1)
  \end{split}
\end{equation}
Let us compute the following term,
\begin{equation}\label{ext-1}
\begin{split}
&df(b_{1, m-1}\otimes e_ja_0)\otimes a_{1, n-1}\otimes f_j\otimes  b_m\otimes 1\\
=&\sum (-1)^{m-1}\delta(df(\id_{m-2}\otimes e_j\otimes a_0)\otimes a_{1, n-1}\otimes f_j\otimes \id\otimes 1)(b_{1, m})+\\
&(-1)^{}\sum df(b_{1, m-2}\otimes e_j\otimes a_0)\otimes
a_{1, n-1}\otimes f_j\otimes b_{m-1}b_m\otimes 1+\\
&\sum df(b_{1, m-2}\otimes e_j\otimes a_0)\otimes
a_{1, n-1}\otimes f_j\otimes b_{m-1}\otimes b_m +\\
&\sum df(b_{1, m-2}\otimes b_{m-1}e_j\otimes a_0)\otimes
a_{1, n-1}\otimes f_j\otimes b_m\otimes 1+\\
&\sum df(b_{1, m-1}\otimes e_j) a_0\otimes
a_{1, n-1}\otimes f_j\otimes b_m\otimes 1\\
=&\sum\sum^{m-1}_{i=1} (-1)^{m-1+(i-1)n} \delta(df(\id_{m-i-1}\otimes e_j\otimes a_{0, i-1})\otimes a_{m-i-1, n-1}\otimes f_j\otimes \id_{m-i}
\otimes 1)(b_{1, m})\\
&\sum \sum_{i=1}^{n} (-1)^{(n-1)i+1}df(b_{1, m-i-1}\otimes e_ja_0\otimes a_{1, i})\otimes a_{i+1
, n-1}\otimes f_j\otimes b_{m-i, m}\otimes 1+\\
&\sum \sum_{i=1}^n(-1)^{i+1}df(b_{1, m-i}\otimes e_j\otimes a_{n-i+1, n-1})a_0\otimes a_{1, n-i}
\otimes f_j\otimes b_{m-i+1, m}\otimes 1+\\
&\sum \sum_{i=1}^{n-1}(-1)^{i(n-1)}e_jf(a_{i, n-1}\otimes a_{0, i-1}\otimes f_j\otimes b_{1, m-n-1})\otimes b_{m-n, m}\otimes 1\\
\end{split}
\end{equation}
Hence the following identity holds by combining (\ref{ext-0}) and (\ref{ext-1}).
\begin{equation}\label{ext3}
\begin{split}
&(\kappa_{n-m+1}(\{f, z\})-[f, \kappa_n(z)])(b_{1, m})\\
=
&(-1)^m \sum de_ja_0\Delta(f)(a_{1, n-1}\otimes f_j\otimes b_{1, m-n-1})\otimes b_{m-n, m}\otimes 1+\\
&(-1)^{m+1}d\Delta(f\cap z)(b_{1, m-n-1})\otimes b_{m-n, m}\otimes 1-\\
     &\sum \sum_{i=1}^{m-n}(-1)^{(m-i)(n+1)}d(f\otimes \id)(b_{1, i-1}\otimes e_ja_0\otimes a_{1, n-1}\otimes f_j\otimes b_{i, m})\otimes 1-\\
   &\sum \sum_{i=1}^{n}(-1)^{(n+1-i)m}d(\id_i\otimes f)(e_ja_0\otimes a_{1, n-1}\otimes f_j\otimes b_{1, m}\otimes 1)-\\
  & \sum \sum_{i=1}^n(-1)^{i+1}df(b_{1, m-i}\otimes e_j\otimes a_{n-i+1, n-1})a_0\otimes a_{1, n-i}
\otimes f_j\otimes b_{m-i+1, m}\otimes 1\\
\end{split}
\end{equation}
By calculation, we have the following identity
\begin{equation}\label{ext4}
  \begin{split}
    &\sum de_ja_0\otimes a_{1, n-1}\otimes f(f_j\otimes b_{1, m-1})\otimes b_m 1\\
    =&\sum de_ja_0\otimes a_{1, n-1}\otimes \langle e_k'f(f_j\otimes b_{1, m-1}), 1\rangle f_k'\otimes b_m \otimes 1\\
    =&\sum\sum_{l=1}^{n}\sum_{i=l}^{m-1} (-1)^{(m+l)(m-i)+1}\delta
    (de_ja_0\otimes a_{1, n-l}\otimes \langle f(\id_{m-i-1}\otimes e_k'\otimes a_{n-l+1, n-1}\otimes f_j\otimes \id_{i-l}), 1\rangle f_k'\\
    &\otimes \id_l\otimes 1)(b_{1, m})+\sum\sum_{i=1}^m(-1)^{m(m-1)} de_ja_0\otimes a_{1, n-1}\otimes
    \Delta(f)(b_{1, m-1})f_j\otimes b_m\otimes 1+\\
    &\sum\sum_{l=1}^{n-1}(-1)^{ml+1}de_ja_0\otimes a_{1,n-l-1}\otimes f(a_{n-l, n-1}
    \otimes f_j\otimes b_{1, m-l-1})\otimes b_{m-l, m}\otimes 1+\\
    &\sum\sum_{i=0}^{n-1}(-1)^{(m+1)i+1}de_ja_0\otimes a_{1, n-i-1}\otimes \langle f(b_{1, m-i-1}\otimes e'_j\otimes a_{n-i, n-1})f_j, 1\rangle f_j'\otimes b_{m-i, m}\otimes 1+\\
    & \sum \sum_{i=1}^{m-n}de_ja_0\langle f(b_{i, m-n-1}\otimes e_j'\otimes a_{1, n-1}\otimes
    f_j\otimes b_{1, i-1}), 1\rangle f_j'\otimes b_{m-n, m}\otimes 1
  \end{split}
\end{equation}
Combining (\ref{ext3}) and (\ref{ext4}), we obtain that
\begin{equation}\label{ext5}
  \begin{split}
    &(\kappa_{n-m+1}(\{f, z\})-[f, \kappa_n(z)])(b_{1, m})\\
=
&(-1)^{m+1}\Delta(f\cap z)(b_{1, m-n-1})\otimes b_{m-n, m}\otimes 1-\\
     &\sum \sum_{i=1}^{m-n}(-1)^{(m-i)(n+1)}d(f\otimes \id)(b_{1, i-1}\otimes e_ja_0\otimes a_{1, n-1}\otimes f_j\otimes b_{i, m})\otimes 1-\\
    & \sum \sum_{i=1}^{m-n}de_ja_0\langle f(b_{i, m-n-1}\otimes e_j'\otimes a_{1, n-1}\otimes
    f_j\otimes b_{1, i-1}), 1\rangle f_j'\otimes b_{m-n, m}\otimes 1
  \end{split}
\end{equation}
Hence from (\ref{ext5}), it remains to verify the following identity,
\begin{equation}\label{ext6}
  \begin{split}
0=&
(-1)^{m+1}\Delta(f\cap z)(b_{1, m-n-1})\otimes b_{m-n, m}\otimes 1-\\
     &\sum \sum_{i=1}^{m-n}(-1)^{(m-i)(n+1)}d(f\otimes \id)(b_{1, i-1}\otimes e_ja_0\otimes a_{1, n-1}\otimes f_j\otimes b_{i, m})\otimes 1-\\
    & \sum \sum_{i=1}^{m-n}de_ja_0\langle f(b_{i, m-n-1}\otimes e_j'\otimes a_{1, n-1}\otimes
    f_j\otimes b_{1, i-1}), 1\rangle f_j'\otimes b_{m-n, m}\otimes 1
  \end{split}
\end{equation}
Let us prove Identity (\ref{ext6}) above.
\begin{equation*}
  \begin{split}
   & \sum d\langle e_ja_0f(a_{1, n-1}\otimes f_j\otimes b_{1, m-n-1}\otimes e_j'), 1\rangle f_j'\otimes b_{m-n, m}\otimes 1\\
    =&\sum(-1)^{n-1} d\langle f(e_ja_0\otimes a_{1, n-1}\otimes f_jb_1\otimes b_{2, m-n-1}\otimes e_j'), 1\rangle f_j'\otimes b_{m-n, m}\otimes 1+\\
    &\sum_{i=1}^{m-n-2} (-1)^{n-1+i}\sum d\langle f(e_ja_0\otimes a_{1, n-1}\otimes f_j\otimes b_{1,i-1} \otimes b_ib_{i+1}\otimes b_{i+2, m-n-1}\otimes e_j'), 1\rangle f_j'\otimes b_{m-n, m}\otimes 1+\\
    &\sum (-1)^{m} d\langle f(e_ja_0\otimes a_{1, n-1}\otimes f_j\otimes b_{1, m-n-2}\otimes b_{m-n-1}e_j'), 1\rangle f_j'\otimes b_{m-n, m}\otimes 1+\\
    &\sum (-1)^{m+1} d\langle f(e_ja_0\otimes a_{1, n-1}\otimes f_j\otimes b_{1, m-n-1})e_j', 1\rangle f_j'\otimes b_{m-n, m}\otimes 1\\
    =&\sum\sum_{i=1}^{m-n}(-1)^{m+i-1}df(b_{1, i}\otimes
    e_ja_0\otimes a_{1, n-1}\otimes f_j\otimes b_{i, m-n-1})\otimes
    b_{m-n, m}\otimes 1+\\
    &\sum \sum^{m-n-1}_{i=1}(-1)^{i(m-n)+1}\langle e_ja_0f(a_{1, n-1}\otimes f_j\otimes b_{i, m-n-1}\otimes e'_j\otimes b_{1, i-1}), 1\rangle f_j'\otimes b_{m-n, m}\otimes 1+\\
    &\sum \sum^{m-n}_{i=1} d\langle f(b_{1, i-1}\otimes e_k'e_ja_0\otimes a_{1, n-1}\otimes f_j\otimes b_{i, m-n-1}), 1\rangle f_k'\otimes b_{m-n,m}\otimes 1.\\
  \end{split}
\end{equation*}
Hence we have the right hand side in (\ref{ext6}) is zero.
\epf

\begin{prop}\label{prop-ext-tor}
Let $A$ be a symmetric algebra over a field $k$. Then we have the following commutative diagram
for $m\geq 2, n\geq 2$,
\begin{equation*}
\xymatrix{
\Ext^1(A, \Omega^{m+1}(A))\otimes \Ext^{1}_{A^e}(A, \Omega^{n+1}(A))\ar[r]^-{[\cdot,\cdot]} &\Ext^1_{A^e}(A, \Omega^{m+n+2}(A))\\
\HH_{\sg}^{-m}(A, A)\otimes \HH_{\sg}^{-n}(A, A)\ar[u]\ar[r]^-{[\cdot, \cdot]}  & \HH_{\sg}^{-m-n-1}(A, A)\ar[u]\\
\Tor_{m-1}^{A^e}(A, A)\otimes \Tor_{n-1}^{A^e}(A, A)\ar[u]^{\kappa_m\otimes \kappa_n}\ar[r]^-{\{\cdot, \cdot\}}  & \Tor^{A^e}_{m+n}(A, A)
\ar[u]_{\kappa_{m+n+1}}
}
\end{equation*}
where $\{\cdot,\cdot\}$ is defined as follows, for any $\alpha \in\Tor_{m-1}^{A^e}(A, A)$ and $\beta\in \Tor_{m-1}^{A^e}(A, A)$
\begin{equation*}
\{\alpha,\beta\}:=(-1)^mB(\alpha)\cup \beta+\alpha\cup B(\beta)+(-1)^{m+1}B(\alpha\cup \beta),
\end{equation*}
where $\cup$ represents the generalized cup product defined in (\ref{gene-cup}).
\end{prop}
\pf The proof is similar to the proof in Proposition \ref{prop-ext1}.
\epf

Therefore, combining Theorem \ref{thm-tra},  Proposition \ref{prop-tor1}, \ref{prop-ext1} and
 \ref{prop-ext-tor}, we obtain the following corollary.
\begin{cor}\label{cor-bv}
  Let $A$ be a symmetric algebra over a field $k$. Then $\HH^*_{\sg}(A, A)$ is a  BV algebra with BV operator $\Delta_{\sg}$,  which is the Connes B-operator for the negative part $\HH^{< 0}_{\sg}(A, A)$, the $\Delta$-operator for the positive part $\HH^{> 0}_{\sg}(A, A)$ and $$\Delta_{\sg}|_{\HH^0_{\sg}(A, A)}=0: \HH^0_{\sg}(A, A)\rightarrow
  \HH^{-1}_{\sg}(A, A).$$ In particular, we have two BV subalgebras $\HH_{\sg}^{\leq 0}(A, A)$ and $\HH_{\sg}^{\geq 0}(A, A)$ with induced BV algebra structures.
\end{cor}
\pf  It remains to prove that we have the following commutative diagram for $m\in \Z_{>0}$,
that is, the image of the bracket $\{\cdot,\cdot\}$ is contained in $\HH_{\sg}^{-1}(A, A)$.
\begin{equation}
\xymatrix{
\HH^m(A, A)\otimes \Tor_{m-1}^{A^e}(A, A)\ar[rd] \ar[r]^-{\{\cdot, \cdot\}} & \Tor_0^{A^e}(A, A)\\
& \HH^{-1}_{\sg}(A, A)\ar@{_(->}[u]
}
\end{equation}
where we recall that the injection
$$\HH^{-1}_{\sg}(A, A)\rightarrow  \Tor_0^{A^e}(A, A)$$
is defined in Proposition \ref{prop-hom} and
$$\{f, \alpha\}:=(-1)^m \Delta(f)\cap \alpha+f\cap B(\alpha)$$
for any $f\in \HH^{m-1}(A, A)$ and $\alpha\in \Tor_{m-1}^{A^e}(A, A)$.
From the short exact sequence in Proposition \ref{prop-hom},
\begin{equation*}
      \xymatrix{
      0\ar[r] & \HH^{-1}_{\sg}(A, A)\ar[r] & A^{\vee}\otimes_{A^e}A\ar[r]^-{\mu^*}&  \Hom_{A^e}(A, A) \ar[r] & \HH^0_{\sg}(A, A)\rightarrow 0
      }
    \end{equation*}
it is sufficient to show that for any $f\in \HH^{m-1}(A, A)$ and $$\alpha:=\sum a_0\otimes a_{1, m-1}\in \Tor_{m-1}^{A^e}(A, A),$$
we have $$\mu^*(\{f, \alpha\})=0.$$
Indeed, we have
\begin{equation*}
\begin{split}
\mu^*(\{f, \alpha\})=&\sum_j (-1)^me_ja_0\Delta(f)(a_{1, m-1}) f_j+\sum_j\sum_{i=0}^{m-1} (-1)^{i(m-1)}e_jf(a_{i, m-1}\otimes a_{0, i-1})f_j.\\
=&\sum_j (-1)^me_ja_0\langle \Delta(f)(a_{1, m-1})f_j e_k', 1\rangle f_k'+ \sum_j\sum_{i=0}^{m-1} (-1)^{i(m-1)}e_jf(a_{i, m-1}\otimes a_{0, i-1})f_j\\
=&\sum_j\sum_{i=1}^m (-1)^{m+i(m-1)}e_ja_0\langle f(a_{i, m-1}\otimes f_je_k'\otimes a_{1, i-1}),  1\rangle f_k'+ \\
&\sum_j\sum_{i=0}^{m-1} (-1)^{i(m-1)}e_jf(a_{i, m-1}\otimes a_{0, i-1})f_j\\
=&0
\end{split}
\end{equation*}
since by direct calculation, we obtain that
\begin{equation*}
\begin{split}
\sum_j\sum_{i=1}^m (-1)^{m+i(m-1)}e_ja_0\langle f(a_{i, m-1}\otimes f_je_k'\otimes a_{1, i-1}),  1\rangle f_k'&=0,\\
\sum_j\sum_{i=0}^{m-1} (-1)^{i(m-1)}e_jf(a_{i, m-1}\otimes a_{0, i-1})f_j&=0.
\end{split}
\end{equation*}
Moreover, we have the following commutative diagram,
\begin{equation*}
\xymatrix{
\HH^m(A, A)\otimes \Tor_{m-1}^{A^e}(A, A) \ar[d]^-{\id\otimes \kappa_m} \ar[r]^-{\{\cdot,\cdot\}} & \HH_{\sg}^{-1}(A, A)\ar@{=}[d]\\
\HH^m(A, A)\otimes \HH^{-m}_{\sg}(A, A) \ar[r]^-{[\cdot, \cdot]} & \HH_{\sg}^{-1}(A, A).
}
\end{equation*}
Therefore, the proof has been completed.
\epf





\begin{cor}\label{cor-cy}
Let $A$ be a symmetric algebra over a field $k$. Then the cyclic homology $\HC_*(A, A)$ is a
graded Lie algebra of lower degree 2, that is, $\HC_*(A, A)[-1]$ is a graded Lie algebra.
\end{cor}
\pf This is an immediate corollary of Porposition 26 in \cite{Men1} since from Corollary \ref{cor-bv}
above it follows that $\HH_*(A, A),$ equipped with Connes B-operator is a BV algebra.
\epf

\bibliographystyle{plain}

\end{document}